\documentclass[8pt]{amsart}

\usepackage{amsmath}
\usepackage{amssymb,latexsym}
\usepackage{amscd}
\usepackage{xy} \xyoption{all} \CompileMatrices
\usepackage[
    bookmarksnumbered = true
    ]{hyperref}
\usepackage{amsthm}
    \newtheorem{theorem}    {Theorem}       [section]
    \newtheorem{lemma}      [theorem]       {Lemma}
    \newtheorem{corollary}  [theorem]       {Corollary}
    \newtheorem{proposition}[theorem]       {Proposition}
    
    \newtheorem{definition} [theorem]       {Definition}
    
    \newtheorem{note}       [theorem]       {Remark}
    
\newcommand{\Z}{{\mathbb Z}}

\newcommand{\R}{{\mathbb R}}
\newcommand{\C}{{\mathbb C}}
\newcommand{\F}{{\mathbb F}}

\newcommand \G {{\mathbb G}}
\newcommand \M {{\mathcal M}}

\newcommand \Spec {{\operatorname{Spec \,}}}
\newcommand \Ord {{\bf{Ord}}}

\usepackage[square] {natbib}

\makeindex

\begin{document}

\title{On multilinearity and skew-symmetry of certain symbols in motivic cohomology of fields}

         \author[Sung Myung]{Sung Myung}
         \address{Department of Mathematics Education, Inha University, 253 Yonghyun-dong, Nam-gu,
             Incheon, 402-751 Korea}
         \email{s-myung1\char`\@inha.ac.kr}

\keywords{motivic cohomology, Milnor's K-theory}
        \date{}
        \maketitle

\begin{abstract}
The purpose of the present article is to show the multilinearity for symbols in Goodwillie-Lichtenbaum complex in two cases.
The first case shown is where the degree is equal to the weight. In this case, the motivic cohomology groups of a field are isomorphic to
the Milnor's $K$-groups as shown by Nesterenko-Suslin, Totaro and Suslin-Voevodsky for various motivic complexes, but we give an explicit isomorphism
for Goodwillie-Lichtenbaum complex in a form which visibly carries multilinearity of Milnor's symbols to our multilinearity of motivic symbols.
Next, we establish multilinearity and skew-symmetry for irreducible Goodwillie-Lichtenbaum symbols in $H^{l-1}_{\M} \bigl(\Spec k , \Z(l) \bigr)$.
These properties have been expected to hold from the author's construction of a bilinear form of dilogarithm
in case $k$ is a subfield of $\C$ and $l=2$.
The multilinearity of symbols may be viewed as a generalization of the well-known formula $\det(AB) = \det(A) \det(B)$ for commuting matrices.
\end{abstract}

\maketitle

\markboth{S. Myung}{multilinearity and skew-symmetry}

       \section{Introduction}
When $R$ is a commutative ring, the group $K_1(R)$ is an abelian group generated by invertible matrices with entries in $R$.
In particular, when $R$ is a field, it is well-known that the determinant map $\det : K_1(R) \rightarrow R^\times$ is an isomorphism.
An important consequence of this fact is that $(AB)=(A)+(B)$, i.e., the product $A B$ of two invertible matrices $A$ and $B$
represents the element obtained by adding two elements in $K_1(R)$,
which are represented by the matrices $A$ and $B$, respectively, since
$\det \begin{pmatrix} A & 0 \\ 0 & B \end{pmatrix} = \det A  \, \det B$. In the present article, we endeavor to generalize this property
to the case of commuting matrices in terms of motivic cohomology. The motivic chain complex proposed
by Goodwillie and Lichtenbaum as follows will be perfectly suitable for our purpose.

In \cite{MR96h:19001}, a chain complex for motivic cohomology of a regular local ring $R$, by Goodwillie and Lichtenbaum,
is defined to be the chain complex associated to the simplicial abelian group $d \mapsto
K_0(R\Delta^d,\, \G_m^{\wedge t})$, together with a shift of degree by $-t$.
Here, $K_0(R\Delta^d,\, \G_m^{\wedge t})$ is the Grothendieck group of the exact category
of projective $R$-modules with $t$ commuting automorphisms factored by the subgroup generated by
classes of the objects one of whose $t$ automorphisms is the identity map.
The motivic cohomology of a regular scheme $X$
is given by hypercohomology of the sheafification of the complex above.
Walker showed, in Theorem 6.5 of \cite{MR2052190}, that it agrees with
motivic cohomology given by Voevodsky and thus various other
definitions of motivic cohomology for smooth schemes over
an algebraically closed field.

In \cite{MR96h:19001}, Grayson showed that a related chain complex
$\Omega^{-t}|d \mapsto K_0^\oplus (R\Delta^d,\, \G_m^{\wedge
t})|$, which uses direct-sum Grothendieck groups instead, arises
as the consecutive quotients in $K$-theory space $K(R)$ when $R$
is a regular noetherian ring and so gives rise to a spectral
sequence converging to $K$-theory. Suslin, in \cite{MR2024054}, showed
that Grayson's motivic cohomology complex is equivalent to the
other definitions of motivic complex and consequently settled the problem of a motivic
spectral sequence. See also \cite{MR2181820} for an overview.

The main results of this article are multilinearity and skew-symmetry properties for the symbols of Goodwillie and Lichtenbaum
in motivic cohomology.
First, we establish them for $H^n_{\M} \bigl(\Spec k , \Z(n) \bigr)$ of a field $k$ in Corollary \ref{multilin-l-l}.
We also give a direct proof of Nesterenko-Suslin's theorem (\cite{MR992981}) that
the motivic cohomology of a field $k$, when the degree is equal to the weight, is equal to the Milnor's $K$-group $K^M_n (k)$
for this version of motivic complex in Theorem \ref{Milnor-iso}.
Even though Nesterenko-Suslin's theorem have already appeared in several articles
including \cite{MR992981}, \cite{MR1187705} and \cite{MR1744945},
we believe that the theorem is a central one in the related subjects and it is worthwhile to have another proof of it.
Moreover, multilinearity and skew-symmetry properties for the symbols of Goodwillie and Lichtenbaum motivic cohomology
$H^n_{\M} \bigl(\Spec k , \Z(n) \bigr)$ and the similar properties for the symbols in Milnor's $K$-groups
are visibly compatible through our isomorphism.
Secondly, we establish multilinearity and skew-symmetry of the irreducible symbols for $H^{l-1}_{\M} \bigl(\Spec k , \Z(l) \bigr)$ in Theorem
\ref{multilinear} and Proposition \ref{skewsymmetry}.
These results are particularly interesting because these are the properties which have been expected
through the construction of the author's regulator map in \cite{MR2189214}
in case $k$ is a subfield of the field $\C$ of complex numbers and $l=2$.
These properties may provide the Goodwillie-Lichtenbaum complex
with a potential to be one of the better descriptions of motivic cohomology of fields.

\section{Multilinearity for Goodwillie-Lichtenbaum motivic complex and Milnor's $K$-groups}

For a ring $R$, let $\mathcal{P}
(R,\, \G_m^l)$ be the exact category each of whose objects
$(P,\theta_1,\dots,\theta_l)$ consists of a finitely generated
projective $R$-module $P$ and commuting automorphisms
$\theta_1,\dots,\theta_l$ of $P$. A morphism from
$(P,\theta_1,\dots,\theta_l)$ to $(P',\theta_1',\dots,\theta_l')$
in this category is a homomorphism $f: P \rightarrow P'$ of
$R$-modules such that $f \theta_i = \theta_i' f$ for each $i$.
Let $K_0(R,\, \G_m^l)$ be the Grothendieck group of this category
and let $K_0(R,\, \G_m^{\wedge l})$ be the quotient of $K_0(R,\,
\G_m^l)$ by the subgroup generated by those objects $(P,\,
\theta_1,\, \dots,\, \theta_l)$ where $\theta_i = 1$ for some $i$.

For each $d \ge 0$, let $R\Delta^d$ be the $R$-algebra
$$R\Delta^d = R[t_0,\dots,t_d]/(t_0 + \cdots + t_d -1).$$
It is isomorphic to a polynomial ring with $d$ indeterminates over
$R$. We denote by \Ord \ the category of finite nonempty ordered
sets and by $[d]$ where $d$ is a nonnegative integer the object
$\{ 0 < 1 < \dots < d\}$. Given a map $\varphi : [d] \rightarrow
[e]$ in \Ord, the map $\varphi^* : R\Delta^e \rightarrow
R\Delta^d$ is defined by $\varphi^*(t_j) = \sum_{\varphi(i)=j}
t_i$. The map $\varphi^*$ gives us a simplicial ring
$R\Delta^\bullet$.

By applying the functor $K_0(-,\, \G_m^{\wedge l})$, we get the simplicial abelian group
$$[d] \mapsto K_0(R\Delta^d,\, \G_m^{\wedge l}).$$
The associated (normalized) chain complex, shifted cohomologically
by $-l$, is called the motivic complex of Goodwillie and
Lichtenbaum of $weight$ $l$.

For each $(P,\theta_1,\dots,\theta_l)$ in $K_0(R,\, \G_m^{\wedge
l})$, there exists a projective module $Q$ such that $P \oplus Q$
is free over $R$. Then $(P \oplus Q,\theta_1 \oplus
1_Q,\dots,\theta_l \oplus 1_Q)$ represents the same element
of $K_0(R,\, \G_m^{\wedge l})$ as $(P,\theta_1,\dots,\theta_l)$.
Thus $K_0(R\Delta^d,\, \G_m^{\wedge l})$ can be explicitly presented
with generators and relations involving $l$-tuples of commuting
matrices in $GL_n(R\Delta^d), \ n \ge 0$.

For a regular local ring $R$, the motivic cohomology $H^q_{\M} \bigl( \Spec R, \, \Z(l) \bigr)$ will be
the $(l-q)$-th homology group of the Goodwillie-Lichtenbaum complex of weight $l$. In particular, when $k$ is any field,
\begin{align*}
H^q_{\M} \bigl( \Spec k, \, \Z(l) \bigr)
= \pi_{l-q} |d \mapsto K_0(k \Delta^d,\, \G_m^{\wedge l})|.
\end{align*}

$K_0 (k \Delta^{d}, \, \G_m^{\wedge l})$ ($l \ge 1$) may be considered as the abelian group
generated by $l$-tuples of the form
$\left( \theta_1(t_1, \dots, t_d),\dots,\theta_l(t_1, \dots, t_d)\right)$
and certain explicit relations, where $\theta_1(t_1, \dots, t_d),\dots,\theta_l(t_1, \dots, t_d)$ are commuting matrices in
$GL_n(k [t_1, \dots, t_d])$ for various $n \ge 1$.

When $d=1$, we set $t=t_1$ and the boundary map $\partial$ on the motivic complex sends
$\left( \theta_1(t),\dots,\theta_l(t) \right)$ in $K_0 (k \Delta^{1}, \, \G_m^{\wedge l})$ to
$\left( \theta_1(1),\dots,\theta_l(1) \right) - \left( \theta_1(0),\dots,\theta_l(0) \right)$ in $K_0 (k \Delta^0, \, \G_m^{\wedge l})$.
We will denote by the same notation $(\theta_1,\dots,\theta_l)$ the element in
$K_0 (k \Delta^0, \, \G_m^{\wedge l}) / \partial K_0 (k \Delta^1, \, \G_m^{\wedge l}) = H^{l}_{\M} \bigl( \Spec k, \, \Z(l) \bigr)$
represented by $(\theta_1,\dots,\theta_l)$, by abuse of notation, whenever $\theta_1,\dots,\theta_l$ are commuting matrices in $GL_n(k)$.

\begin{lemma} \label{basicelements}
Let $a_1,a_2, \dots, a_n$ and $b_1, b_2,\dots, b_n$ be elements in $\bar k$ (an algebraic closure of $k$)
not equal to either 0 or 1. Suppose also that $a_1 a_2 \cdots a_n = b_1 b_2 \cdots b_n$ and
$(1-a_1)(1-a_2) \cdots (1-a_n) = (1-b_1)(1-b_2) \cdots (1-b_n)$.
If all the elementary symmetric functions evaluated at $a_1, a_2, \dots , a_n$ and $b_1, b_2, \dots , b_n$
are in $k$, then there is a matrix $\theta(t)$
in $GL_n(k[t])$ such that $1_n-\theta(t)$ is also invertible and the eigenvalues of $\theta(0)$
and $\theta(1)$ are $a_1,a_2, \dots, a_n$ and $b_1,b_2, \dots, b_n$, respectively.
\end{lemma}

\begin{proof}
Let
$$p(\lambda) = (1-t) \prod_{i=1}^n (\lambda -a_i) + t \prod_{i=1}^n (\lambda -b_i)$$
be a polynomial in $\lambda$ with coefficients in $k[t]$.
It is a monic polynomial with the constant term equal to $(-1)^n a_1 a_2 \cdots a_n$.
It has roots $b_1, b_2, \dots, b_n$ and $a_1,a_2, \dots, a_n$ when $t=1$ and $t=0$, respectively.

Now let $\theta(t)$ be its companion matrix in $GL_n(k[t])$. Then $\text{det} \, (1_n-\theta(t)) = p(1)$ since
$\text{det} \, (\lambda 1_n-\theta(t)) = p(\lambda)$.
But $p(1)=(1-a_1)(1-a_2) \cdots (1-a_n)$ $ = (1-b_1)(1-b_2) \cdots (1-b_n)$ is in $k^{\times}$,
and so $1_n-\theta(t)$ is invertible. It is clear that the eigenvalues of $\theta(t)$ are $a_1,\, a_2, \,  \dots, \, a_n$
and $b_1,\, b_2, \,  \dots, \, b_n$ when $t=0$ and $t=1$, respectively.
\end{proof}

\begin{definition} \label{Z}

For $l \ge 2$, let $\mathbf{Z}$ be the subgroup of $K_0 (k \Delta^1, \G_m^{\wedge l})$ generated by the elements of the following types
for various $n \ge 1$ :

$(Z_1)$ $(\theta_1,\dots,\theta_l)$, where $\theta_1,\dots,\theta_l \in GL_n(k[t])$ commute and $\theta_i$ is in $GL_n(k)$ for some $i$;

$(Z_2)$ $(\theta_1,\dots,\theta_l)$, where $\theta_i = \theta_j \in GL_n(k[t])$ for some $i \ne j$;

$(Z_3)$ $(\theta_1,\dots,\theta_l)$, where $\theta_i = 1_n -\theta_j \in GL_n(k[t])$ for some $i \ne j$.

\end{definition}

\begin{lemma} \label{boundary-Z}
Let $\partial \mathbf{Z}$ denote the image of $\mathbf{Z}$ under the boundary homomorphism
$\partial: K_0 (k \Delta^1, \, \G_m^{\wedge l}) \rightarrow K_0 (k \Delta^0, \, \G_m^{\wedge l})$ when $l \ge 2$.
Then $\partial \mathbf{Z}$ contains all elements of the following forms:

(i) $(\varphi \psi,\theta_2,\dots,\theta_l)-(\varphi,\theta_2,\dots,\theta_l)-(\psi,\theta_2,\dots,\theta_l)$,
for all commuting $\varphi, \psi, \theta_2,\dots,\theta_l \in GL_n(k)$;

 Similarly, $(\theta_1, \dots, \theta_{i-1}, \varphi \psi, \theta_{i+1}, \dots, \theta_l)
 -(\theta_1, \dots, \theta_{i-1}, \varphi, \theta_{i+1}, \dots, \theta_l)
 -(\theta_1, \dots, \theta_{i-1}, \psi, \theta_{i+1}, \dots, \theta_l)$
for all commuting $\varphi, \psi, \theta_1,\dots,\theta_{i-1}, \theta_{i+1},\dots, \theta_l \in GL_n(k)$;

(ii) $(\theta_1,\dots,\theta_i, \dots, \theta_j, \dots, \theta_l) + (\theta_1,\dots,\theta_j, \dots, \theta_i, \dots, \theta_l)$,
for all commuting $\theta_1,\dots,\theta_l \in GL_n(k)$;

(iii) $(\theta_1,\dots,\theta_i, \dots, \theta_j, \dots, \theta_l)$, when $\theta_i= -\theta_j$
for commuting $\theta_1,\dots,\theta_l \in GL_n(k)$;

(iv) $(c_1, \dots, b,\dots, 1-b, \dots, c_l)-(c_1, \dots, a,\dots, 1-a,\dots, c_l)$, for $a, b \in k-\{0,1\}$ and $c_i \in k^\times$
for each appropriate $i$.
\end{lemma}

\begin{proof}
$(i)$
We first observe the following identities of matrices:
\begin{align} \label{linear-equa}
\begin{pmatrix} 1_n & 0 \\ \psi & 1_n \end{pmatrix}  \begin{pmatrix} \psi & 1_n \\ 0 & \varphi \end{pmatrix}  \begin{pmatrix} 1_n & 0 \\ -\psi & 1_n \end{pmatrix}
&=\begin{pmatrix} 0 & 1_n \\ -\varphi \psi & \varphi+\psi \end{pmatrix}, \\
\label{linear-equb}
\begin{pmatrix} 1_n & 0 \\ 1_n & 1_n \end{pmatrix}
 \begin{pmatrix} 1_n & 1_n \\ 0 & \varphi \psi \end{pmatrix}
 \begin{pmatrix} 1_n & 0 \\ -1_n & 1_n \end{pmatrix}
&=\begin{pmatrix} 0 & 1_n \\ -\varphi \psi & 1_n+\varphi \psi \end{pmatrix}.
\end{align}

Let $\Theta(t)$ be the $2n \times 2n$ matrix
$$\begin{pmatrix} 0 & 1_n \\
              -\varphi \psi & t(1_n+\varphi \psi)+(1-t)(\varphi+\psi) \end{pmatrix}. $$
Then, $\Theta(t)$ is in $GL_{2n}(k[t])$, $\bigl( \Theta(t),\, \theta_2 \oplus \theta_2, \dots, \theta_l \oplus \theta_l \bigr)$
is in $\mathbf{Z}$ by Definition \ref{Z} $(Z_1)$
and the boundary of $\bigl( \Theta(t),\, \theta_2 \oplus \theta_2, \dots, \theta_l \oplus \theta_l \bigr)$ is, by (\ref{linear-equa})
and by (\ref{linear-equb}),
$(1_n \oplus \varphi \psi, \theta_2 \oplus \theta_2, \dots, \theta_l \oplus \theta_l)
-(\varphi \oplus \psi,\theta_2 \oplus \theta_2, \dots, \theta_l \oplus \theta_l)
=(\varphi \psi,\theta_2, \dots, \theta_l)
-(\varphi,\theta_2, \dots, \theta_l)-(\psi,\theta_2, \dots, \theta_l)$.

The proof is similar for other cases.

$(ii)$ We let $\Theta(t)$ be the matrix
$$\begin{pmatrix} 0 & 1_n \\
              -\theta_i \theta_j & t(1_n+\theta_i \theta_j)+(1-t)(\theta_i+\theta_j) \end{pmatrix}. $$

 Then $\bigl( \theta_1^{\oplus 2},\dots,\Theta(t), \dots, \Theta(t),\dots, \theta_l^{\oplus 2} \bigr)$ is in $\mathbf{Z}$
by Definition \ref{Z} $(Z_2)$ and the boundary of $\bigl( \theta_1^{\oplus 2},\dots,\Theta(t), \dots, \Theta(t),\dots, \theta_l^{\oplus 2} \bigr)$ is
\begin{align*}
&(\theta_1,\dots,\theta_i \theta_j, \dots, \theta_i \theta_j,\dots, \dots, \theta_l)
-(\theta_1,\dots,\theta_i, \dots, \theta_i,\dots, \theta_l)
-(\theta_1,\dots,\theta_j, \dots, \theta_j,\dots,  \theta_l)\\
&= \bigl( (\theta_1,\dots,\theta_i, \dots, \theta_i,\dots, \theta_l)
+((\theta_1,\dots,\theta_i, \dots, \theta_j,\dots, \theta_l)
+(\theta_1,\dots,\theta_j, \dots, \theta_i,\dots, \theta_l)
+(\theta_1,\dots,\theta_j, \dots, \theta_j,\dots, \theta_l) \bigr) \\
&-(\theta_1,\dots,\theta_i, \dots, \theta_i,\dots,  \theta_l)
-(\theta_1,\dots,\theta_j, \dots, \theta_j,\dots, \theta_l) \\
&=(\theta_1,\dots,\theta_j, \dots, \theta_i,\dots, \theta_l) + (\theta_1,\dots,\theta_j, \dots, \theta_i,\dots, \theta_l)
\quad \text {modulo } \partial \mathbf{Z} \text{ by} \ (i).
\end{align*}

$(iii)$ We note that $\displaystyle \left( \begin{pmatrix} \theta_1 & 0 \\ 0 & \theta_1 \end{pmatrix}, \dots,
\begin{pmatrix} -\theta & 0 \\ 0 & -\theta \end{pmatrix},\dots, \begin{pmatrix} 0 & 1_n \\ -\theta & t(\theta+1_n) \end{pmatrix}
        \dots,  \begin{pmatrix} \theta_l & 0 \\ 0 & \theta_l \end{pmatrix}  \right)$
is an element of $\mathbf{Z}$ by Definition \ref{Z} $(Z_1)$. So its boundary

\begin{multline*}
\quad \ \left( \begin{pmatrix} \theta_1 & 0 \\ 0 & \theta_1 \end{pmatrix}, \dots,
\begin{pmatrix} -\theta & 0 \\ 0 & -\theta \end{pmatrix}, \dots, \begin{pmatrix} 0 & 1_n \\ -\theta & \theta+1_n \end{pmatrix},
           \dots,  \begin{pmatrix} \theta_l & 0 \\ 0 & \theta_l \end{pmatrix} \right) \\
-\left( \begin{pmatrix} \theta_1 & 0 \\ 0 & \theta_1 \end{pmatrix}, \dots,
\begin{pmatrix} -\theta & 0 \\ 0 & -\theta \end{pmatrix}, \dots, \begin{pmatrix} 0 & 1_n \\ -\theta & 0 \end{pmatrix}
 ,\dots,  \begin{pmatrix} \theta_l & 0 \\ 0 & \theta_l \end{pmatrix} \right) \\
=\left(\begin{pmatrix} \theta_1 & 0 \\ 0 & \theta_1 \end{pmatrix}, \dots,
\begin{pmatrix} -\theta & 0 \\ 0 & -\theta \end{pmatrix}, \dots, \begin{pmatrix} \theta & 1_n \\ 0 & 1_n \end{pmatrix}
 ,\dots,  \begin{pmatrix} \theta_l & 0 \\ 0 & \theta_l \end{pmatrix} \right)\\
-\left(\begin{pmatrix} \theta_1 & 0 \\ 0 & \theta_1 \end{pmatrix}, \dots,
\begin{pmatrix} -\theta & 0 \\ 0 & -\theta \end{pmatrix}, \dots, \begin{pmatrix} 0 & 1_n \\ -\theta & 0 \end{pmatrix}
 ,\dots,  \begin{pmatrix} \theta_l & 0 \\ 0 & \theta_l \end{pmatrix} \right) \\
=(\theta_1, \dots, -\theta, \dots, ,\theta, \dots, \theta_l)
-\left(\begin{pmatrix} \theta_1 & 0 \\ 0 & \theta_1 \end{pmatrix}, \dots,
\begin{pmatrix} -\theta & 0 \\ 0 & -\theta \end{pmatrix}, \dots, \begin{pmatrix} 0 & 1_n \\ -\theta & 0 \end{pmatrix}
 ,\dots,  \begin{pmatrix} \theta_l & 0 \\ 0 & \theta_l \end{pmatrix}\right)
\end{multline*}
is in $\partial \mathbf{Z}$. Thus it suffices to prove that
$\displaystyle \left(\begin{pmatrix} \theta_1 & 0 \\ 0 & \theta_1 \end{pmatrix}, \dots,
\begin{pmatrix} -\theta & 0 \\ 0 & -\theta \end{pmatrix}, \dots, \begin{pmatrix} 0 & 1_n \\ -\theta & 0 \end{pmatrix}
 ,\dots,  \begin{pmatrix} \theta_l & 0 \\ 0 & \theta_l \end{pmatrix} \right)$
is in $\partial \mathbf{Z}$. But it is equal to
\begin{multline*}
\left(\begin{pmatrix} \theta_1 & 0 \\ 0 & \theta_1 \end{pmatrix}, \dots,
{\begin{pmatrix} 0 & 1_n \\ -\theta & 0 \end{pmatrix}}^2, \dots, \begin{pmatrix} 0 & 1_n \\ -\theta & 0 \end{pmatrix}
 ,\dots,  \begin{pmatrix} \theta_l & 0 \\ 0 & \theta_l \end{pmatrix} \right) \\
= 2 \left(\begin{pmatrix} \theta_1 & 0 \\ 0 & \theta_1 \end{pmatrix}, \dots,
{\begin{pmatrix} 0 & 1_n \\ -\theta & 0 \end{pmatrix}}, \dots, \begin{pmatrix} 0 & 1_n \\ -\theta & 0 \end{pmatrix}
 ,\dots,  \begin{pmatrix} \theta_l & 0 \\ 0 & \theta_l \end{pmatrix} \right),
\end{multline*}
which is in $\partial \mathbf{Z}$ by $(ii)$ above.

$(iv)$ Apply Lemma \ref{basicelements} to $a_1=a,\ a_2=\sqrt{b},\ a_3=-\sqrt{b},
 \ b_1=-\sqrt{a}, \ b_2=\sqrt{a}, \ b_3= b$ to get $\theta(t) \in GL_3(k[t])$ with the properties stated in the lemma.
Then $z=2 \bigl( c_1^{\oplus 3}, \dots, \theta(t),\dots, 1_3-\theta(t),\dots,c_l^{\oplus 3} \bigr)$ is in $\mathbf{Z}$ by Definition \ref{Z} $(Z_3)$. But, by the theory of
rational canonical form, we have
{\allowdisplaybreaks
\begin{multline*}
\partial z = 2 \left( (c_1, \dots,b,\dots,1-b,\dots,c_l)
+\left( \begin{pmatrix} c_1 & 0 \\ 0 & c_1 \end{pmatrix}, \dots,
\begin{pmatrix} 0 & 1 \\ a & 0 \end{pmatrix}, \dots, \begin{pmatrix} 1 & -1 \\ -a & 1 \end{pmatrix}
,\dots, \begin{pmatrix} c_l & 0 \\ 0 & c_l \end{pmatrix}\right) \right) \\
-2 \left( (c_1, \dots,a,\dots,1-a,\dots,c_l)
+\left(\begin{pmatrix} c_1 & 0 \\ 0 & c_1 \end{pmatrix}, \dots,
\begin{pmatrix} 0 & 1 \\ b & 0 \end{pmatrix}, \dots, \begin{pmatrix} 1 & -1 \\ -b & 1 \end{pmatrix}
,\dots, \begin{pmatrix} c_l & 0 \\ 0 & c_l \end{pmatrix}\right) \right) \\
= -2 (c_1, \dots,a,\dots,1-a,\dots,c_l) + 2(c_1, \dots,b,\dots,1-b,\dots,c_l) \\
-\left(\begin{pmatrix} c_1 & 0 \\ 0 & c_1 \end{pmatrix}, \dots,
{\begin{pmatrix} 0 & 1 \\ b & 0 \end{pmatrix}}^2, \dots, \begin{pmatrix} 1 & -1 \\ -b & 1 \end{pmatrix}
,\dots, \begin{pmatrix} c_l & 0 \\ 0 & c_l \end{pmatrix}\right) \\
+\left(\begin{pmatrix} c_1 & 0 \\ 0 & c_1 \end{pmatrix}, \dots,
{\begin{pmatrix} 0 & 1 \\ a & 0 \end{pmatrix}}^2, \dots, \begin{pmatrix} 1 & -1 \\ -a & 1 \end{pmatrix}
,\dots, \begin{pmatrix} c_l & 0 \\ 0 & c_l \end{pmatrix}\right)\\
= \left(\begin{pmatrix} c_1 & 0 \\ 0 & c_1 \end{pmatrix}, \dots,
\begin{pmatrix} b & 0 \\ 0 & b \end{pmatrix}, \dots, \begin{pmatrix} 1-b & 0 \\ 0 & 1-b \end{pmatrix}
,\dots, \begin{pmatrix} c_l & 0 \\ 0 & c_l \end{pmatrix}\right) \\
- \left(\begin{pmatrix} c_1 & 0 \\ 0 & c_1 \end{pmatrix}, \dots,
 \begin{pmatrix} b & 0 \\ 0 & b \end{pmatrix}, \dots, \begin{pmatrix} 1 & -1 \\ -b & 1 \end{pmatrix}
,\dots, \begin{pmatrix} c_l & 0 \\ 0 & c_l \end{pmatrix} \right) \\
 -\left(\begin{pmatrix} c_1 & 0 \\ 0 & c_1 \end{pmatrix}, \dots,
\begin{pmatrix} a & 0 \\ 0 & a \end{pmatrix}, \dots, \begin{pmatrix} 1-a & 0 \\ 0 & 1-a \end{pmatrix}
,\dots, \begin{pmatrix} c_l & 0 \\ 0 & c_l \end{pmatrix}\right) \\
 + \left(\begin{pmatrix} c_1 & 0 \\ 0 & c_1 \end{pmatrix}, \dots,
\begin{pmatrix} a & 0 \\ 0 & a \end{pmatrix}, \dots, \begin{pmatrix} 1 & -1 \\ -a & 1 \end{pmatrix}
,\dots, \begin{pmatrix} c_l & 0 \\ 0 & c_l \end{pmatrix}\right) \\
=\left(\begin{pmatrix} c_1 & 0 \\ 0 & c_1 \end{pmatrix}, \dots,
\begin{pmatrix} b & 0 \\ 0 & b \end{pmatrix}, \dots,
\begin{pmatrix} 1-b & 0 \\ 0 & 1-b \end{pmatrix} {\begin{pmatrix} 1 & -1 \\ -b & 1 \end{pmatrix}}^{-1}
  ,\dots, \begin{pmatrix} c_l & 0 \\ 0 & c_l \end{pmatrix} \right) \\
 -\left(\begin{pmatrix} c_1 & 0 \\ 0 & c_1 \end{pmatrix}, \dots,
\begin{pmatrix} a & 0 \\ 0 & a \end{pmatrix}, \dots,
 \begin{pmatrix} 1-a & 0 \\ 0 & 1-a \end{pmatrix} {\begin{pmatrix} 1 & -1 \\ -a & 1 \end{pmatrix}}^{-1}
  ,\dots, \begin{pmatrix} c_l & 0 \\ 0 & c_l \end{pmatrix}  \right)\\
=\left(\begin{pmatrix} c_1 & 0 \\ 0 & c_1 \end{pmatrix}, \dots,
\begin{pmatrix} b & 0 \\ 0 & b \end{pmatrix}, \dots \begin{pmatrix} 1 & 1 \\ b & 1 \end{pmatrix}
  ,\dots, \begin{pmatrix} c_l & 0 \\ 0 & c_l \end{pmatrix} \right)
    -\left(\begin{pmatrix} c_1 & 0 \\ 0 & c_1 \end{pmatrix}, \dots,
\begin{pmatrix} a & 0 \\ 0 & a \end{pmatrix}, \dots, \begin{pmatrix} 1 & 1 \\ a & 1 \end{pmatrix}
 ,\dots, \begin{pmatrix} c_l & 0 \\ 0 & c_l \end{pmatrix} \right) \\
=\left(\begin{pmatrix} c_1 & 0 \\ 0 & c_1 \end{pmatrix}, \dots,
\begin{pmatrix} b & 0 \\ 0 & b \end{pmatrix},\dots, \begin{pmatrix} {\frac {-b} {1-b}} & {\frac 1 {1-b}} \\ 0 & 1 \end{pmatrix}
            \begin{pmatrix} 1 & 1 \\ b & 1 \end{pmatrix} {\begin{pmatrix} {\frac {-b}{1-b}} & {\frac 1 {1-b}} \\ 0 & 1 \end{pmatrix}}^{-1}
          ,\dots, \begin{pmatrix} c_l & 0 \\ 0 & c_l \end{pmatrix}  \right)\\
   -\left(\begin{pmatrix} c_1 & 0 \\ 0 & c_1 \end{pmatrix}, \dots,
\begin{pmatrix} a & 0 \\ 0 & a \end{pmatrix},\dots, \begin{pmatrix} {\frac {-a} {1-a}} & {\frac 1 {1-a}} \\ 0 & 1 \end{pmatrix}
           \begin{pmatrix} 1 & 1 \\ a & 1 \end{pmatrix} {\begin{pmatrix} {\frac {-a}{1-a}} & {\frac 1 {1-a}} \\ 0 & 1 \end{pmatrix}}^{-1}
          ,\dots, \begin{pmatrix} c_l & 0 \\ 0 & c_l \end{pmatrix} \right)\\
=\left(\begin{pmatrix} c_1 & 0 \\ 0 & c_1 \end{pmatrix}, \dots,
\begin{pmatrix} b & 0 \\ 0 & b \end{pmatrix}, \dots, \begin{pmatrix} 0 & 1 \\ b-1 & 2 \end{pmatrix}
,\dots, \begin{pmatrix} c_l & 0 \\ 0 & c_l \end{pmatrix} \right)\\
-\left(\begin{pmatrix} c_1 & 0 \\ 0 & c_1 \end{pmatrix}, \dots,
\begin{pmatrix} a & 0 \\ 0 & a \end{pmatrix}, \dots, \begin{pmatrix} 0 & 1 \\ a-1 & 2 \end{pmatrix}
 ,\dots, \begin{pmatrix} c_l & 0 \\ 0 & c_l \end{pmatrix}   \right).
\end{multline*}
}
By taking the boundary of the element
\begin{align*}
&\left(\begin{pmatrix} c_1 & 0 \\ 0 & c_1 \end{pmatrix}, \dots,
 \begin{pmatrix} b & 0 \\ 0 & b \end{pmatrix}, \dots, \begin{pmatrix} 0 & 1 \\ b-1 & (2-b)t+2(1-t) \end{pmatrix},
  \dots, \begin{pmatrix} c_l & 0 \\ 0 & c_l \end{pmatrix} \right) \\
&-\left(\begin{pmatrix} c_1 & 0 \\ 0 & c_1 \end{pmatrix}, \dots,
\begin{pmatrix} a & 0 \\ 0 & a \end{pmatrix}, \dots, \begin{pmatrix} 0 & 1 \\ a-1 & (2-a)t+2(1-t) \end{pmatrix},
 \dots, \begin{pmatrix} c_l & 0 \\ 0 & c_l \end{pmatrix} \right),
\end{align*}
which is in $\mathbf{Z}$ by Definition \ref{Z} $(Z_1)$, we see that
\begin{align*}
\partial z
&= \left(\begin{pmatrix} c_1 & 0 \\ 0 & c_1 \end{pmatrix}, \dots,
\begin{pmatrix} b & 0 \\ 0 & b \end{pmatrix}, \dots, \begin{pmatrix} 0 & 1 \\ b-1 & 2-b \end{pmatrix},
 \dots, \begin{pmatrix} c_l & 0 \\ 0 & c_l \end{pmatrix} \right)\\
&-\left(\begin{pmatrix} c_1 & 0 \\ 0 & c_1 \end{pmatrix}, \dots,
\begin{pmatrix} a & 0 \\ 0 & a \end{pmatrix}, \dots, \begin{pmatrix} 0 & 1 \\ a-1 & 2-a \end{pmatrix},
 \dots, \begin{pmatrix} c_l & 0 \\ 0 & c_l \end{pmatrix}\right) \\
&= \left(\begin{pmatrix} c_1 & 0 \\ 0 & c_1 \end{pmatrix}, \dots,
\begin{pmatrix} b & 0 \\ 0 & b \end{pmatrix}, \dots, \begin{pmatrix} 1-b & 0 \\ 0 & 1 \end{pmatrix},
 \dots, \begin{pmatrix} c_l & 0 \\ 0 & c_l \end{pmatrix}\right)
\end{align*}
by (\ref{linear-equa}), which then is equal to
\begin{align*}
& -\left(\begin{pmatrix} c_1 & 0 \\ 0 & c_1 \end{pmatrix}, \dots,
    \begin{pmatrix} a & 0 \\ 0 & a \end{pmatrix}, \dots, \begin{pmatrix} 1-a & 0 \\ 0 & 1 \end{pmatrix},
    \dots, \begin{pmatrix} c_l & 0 \\ 0 & c_l \end{pmatrix}  \right)\\
&= (c_1, \dots, b,\dots, 1-b, \dots, c_l)-(c_1, \dots, a,\dots, 1-a,\dots, c_l) \\
\end{align*}
in $K_0 (k \Delta^0, \, \G_m^{\wedge l})/\partial \mathbf{Z}$.
Therefore, $(iv)$ lies in $\partial \mathbf{Z}$.
\end{proof}

\begin{corollary} \label{multilin-l-l} (Multilinearity and Skew-symmetry for $H^l_{\M} \bigl( \Spec k, \, \Z(l) \bigr)$)

(i) $(\theta_1, \dots, \theta_{i-1}, \varphi \psi, \theta_{i+1}, \dots, \theta_l) =
 (\theta_1, \dots, \theta_{i-1}, \varphi, \theta_{i+1}, \dots, \theta_l) +(\theta_1, \dots, \theta_{i-1}, \psi, \theta_{i+1}, \dots, \theta_l)$
 in $H^l_{\M} \bigl( \Spec k, \, \Z(l) \bigr)$,
for all commuting $\varphi, \psi, \theta_1,\dots,\theta_{i-1}, \theta_{i+1},\dots, \theta_l \in GL_n(k)$

(ii) $(\theta_1,\dots,\theta_i, \dots, \theta_j, \dots, \theta_l) = - (\theta_1,\dots,\theta_j, \dots, \theta_i, \dots, \theta_l)$
in $H^l_{\M} \bigl( \Spec k, \, \Z(l) \bigr)$
for all commuting $\theta_1,\dots,\theta_l \in GL_n(k)$
\end{corollary}

If $\theta_1,\dots,\theta_l$ and $\theta'_1,\dots,\theta'_l$ are commuting matrices in $GL_n(k)$ and $GL_m(k)$, respectively, then
$(\theta_1,\dots,\theta_l) + (\theta'_1,\dots,\theta'_l) = (\theta_1 \oplus \theta'_1,\dots,\theta_l \oplus \theta'_l)$
in $H^l_{\M} \bigl( \Spec k, \, \Z(l) \bigr)$. Therefore, we obtain the following result from Corollary \ref{multilin-l-l}.
\begin{corollary}
Every element in $H^l_{\M} \bigl( \Spec k, \, \Z(l) \bigr)$ can be written as a single symbol $(\theta_1,\dots,\theta_l)$, where
$\theta_1,\dots,\theta_l$ are commuting matrices in $GL_n(k)$.
\end{corollary}

Thanks to Lemma \ref{boundary-Z}, we can construct a map from Milnor's $K$-groups to the motivic cohomology groups.

\begin{proposition} \label{Milnor-map}
For any field $k$, the assignment $\{ a_1,a_2, \dots, a_l \} \mapsto (a_1,a_2,\dots,a_l) $ for each Steinberg symbol
$\{a_1,a_2, \dots, a_l \}$ gives a well-defined homomorphism $\rho_l$
from the Milnor's $K$-group $K^M_l(k)$ to $H^l_{\M} \bigl( \Spec k, \, \Z(l) \bigr)$.
\end{proposition}
\begin{proof}
This proposition turns out to be straightforward when $l=1$. So we assume that $l \ge 2$.
By Corollary \ref{multilin-l-l} $(i)$, the multilinearity is satisfied by our symbol $( \ ,\dots, \ )$.
Therefore all we need to show is that for every $\alpha \in k-\{0,1\}$ and $c_r \in k^\times$ for $1 \le r \le l$, $r \ne i,j$,
$(c_1, \dots, \alpha, \dots, 1-\alpha, \dots, c_l)$ is in $\partial K_0 (k \Delta^1, \, \G_m^{\wedge l})$.
We will actually show that it is contained in $\partial \mathbf{Z}$.

The proposition is immediate for a prime field $\F_p$ because $K^M_l(\F_p)=0$ for $l \ge 2$.
So we may assume that there exists an element $e \in k$ such that $e^3-e \neq 0$.
By Lemma \ref{boundary-Z} $(iv)$ with $a=e,\, b=1-e$, we have
$(c_1, \dots, e,\dots,1-e, \dots, c_l)-(c_1, \dots, 1-e,\dots,e, \dots, c_l) = 2 (c_1, \dots, e,\dots,1-e, \dots, c_l) = 0$
modulo $\partial \mathbf{Z}$.
With $a=-e,\, b=1+e$, we have $2(c_1, \dots, e,\dots,1+e, \dots, c_l) = 2(c_1, \dots,-e,\dots,1+e, \dots, c_l)=0$.
Hence, $(c_1, \dots, e^2,\dots,1-e^2, \dots, c_l) = 2(c_1, \dots, e,\dots,1-e, \dots, c_l)+2(c_1, \dots, e,\dots,1+e, \dots, c_l) = 0$.

On the other hand, by Lemma \ref{boundary-Z} $(iv)$ with $a=e^2,\, b=\alpha$,
we see that $-(c_1, \dots, e^2,\dots,1-e^2, \dots, c_l) + (c_1, \dots, \alpha,\dots,1-\alpha, \dots, c_l)$ is in
$\partial \mathbf{Z}$ and we're done.

More explicitly, let $z = 2 \bigl(c_1^{\oplus 3}, \dots, \theta(t), \dots, 1-\theta(t), \dots, c_l^{\oplus 3} \bigr) \in \mathbf{Z}$, where
$$\theta(t) = \begin{pmatrix} 0 & 1 & 0 \\
            0 &  0  & 1 \\
            -e^2\alpha  &  (e^2-\alpha)t +\alpha & (\alpha-e^2)t+e^2 \end{pmatrix}.$$
This matrix $\theta(t)$ is constructed with Lemma \ref{basicelements} with $a_1=e^2,\ a_2=\sqrt{\alpha},\ a_3=-\sqrt{\alpha},
 \ b_1=-e, \ b_2=e, \ b_3= \alpha$. Hence, by the computation we have done in the proof of Lemma \ref{boundary-Z} $(iv)$,
\begin{align*}
\partial z
&= 2(c_1, \dots, -e,\dots,1+e, \dots, c_l)+2(c_1, \dots, e,\dots,1-e, \dots, c_l) \\
&+ 2(c_1, \dots, \alpha,\dots,1-\alpha, \dots, c_l)-2(c_1, \dots, e^2,\dots,1-e^2, \dots, c_l) \\
&-2 \left(\begin{pmatrix} c_1 & 0 \\ 0 & c_1 \end{pmatrix}, \dots,
\begin{pmatrix} 0 & 1 \\ \alpha & 0 \end{pmatrix}, \dots, \begin{pmatrix} 1 & -1 \\ -\alpha & 1 \end{pmatrix}
 \dots, \begin{pmatrix} c_l & 0 \\ 0 & c_l \end{pmatrix} \right) \\
&=-(c_1, \dots, e^2,\dots,1-e^2, \dots, c_l) + (c_1, \dots, \alpha,\dots,1-\alpha, \dots, c_l) \\
&=((c_1, \dots, \alpha,\dots,1-\alpha, \dots, c_l).
\end{align*}
\end{proof}

For Goodwillie-Lichtenbaum motivic complex, there is a straightforward \underline{functorial} definition of the norm map for the motivic cohomology
for any finite extension $k \subset L$.

\begin{definition} \label{def-norm}
If $\theta_1, \dots, \theta_l$ are commuting automorphisms on a finitely generated projective $L \Delta^d$-module $P$,
then by identifying $L \Delta^d$ as a free $k \Delta^d$-module of finite rank,
we may consider $P$ as a finitely generated projective $k \Delta^d$-module and $\theta_1, \dots, \theta_l$ as commuting automorphisms on it.
This gives a simplicial map $K_0(L \Delta^d,\, \G_m^{\wedge l}) \rightarrow K_0(k \Delta^d,\, \G_m^{\wedge l})$.
The resulting homomorphism $N_{L/k}: \ H^{q}_{\M} \bigl(\Spec L , \Z(l) \bigr) \rightarrow H^{q}_{\M} \bigl(\Spec k , \Z(l) \bigr)$ is called the
norm map.
\end{definition}

We summarize some basic results for the norm in the following lemma.

\begin{lemma} \label{basic-norm}
(i) $N_{L'/L} \circ N_{L/k} = N_{L'/k}$ whenever we have a tower of finite field extensions $k \subset L \subset L'$.

(ii) If $[L:k]=d$, the composition
$$ \xymatrix {H^{q}_{\M} \bigl(\Spec k , \Z(l) \bigr) \ar[r]^-{i_{L/k}} & H^{q}_{\M} \bigl(\Spec L , \Z(l) \bigr) \ar[r]^-{N_{L/k}}
& H^{q}_{\M} \bigl(\Spec k , \Z(l) \bigr) },$$ where $i_{L/k}$ is induced by the inclusion of the fields $k \subset L$,
is multiplication by $d$.

(iii) For $\alpha_1, \dots, \alpha_l \in k^\times$ and $\beta \in L^\times$,
$N_{L/k} \left( \alpha_1, \dots, \alpha_l, \beta \right) = \left(\alpha_1, \dots, \alpha_l, N_{L/k}(\beta) \right)$
in $H^{l+1}_{\M} \bigl(\Spec k , \Z(l+1) \bigr)$,
where $N_{L/k}(\beta) \in k^\times$ is the image of $\beta$ under the usual norm map $N_{L/k}: L^\times \rightarrow k^\times$.
\end{lemma}
\begin{proof}
$(i)$ and $(ii)$ are immediate from Definition \ref{def-norm}. $(iii)$ follows from the observation that, in $H^{1}_{\M} \bigl(\Spec k , \Z(1) \bigr)$,
the two elements represented by two matrices with same determinants are equal since any matrix with determinant 1 is a product
of elementary matrices and an element represented by an elementary matrix vanishes in $H^{1}_{\M} \bigl(\Spec k , \Z(1) \bigr)$.
\end{proof}

We also have the norm maps $N_{L/K}: \ K^M_l (L) \rightarrow K^M_l (k)$ for the Milnor's $K$-groups
whenever $L/k$ is a finite field extension, whose definition we recall briefly as follows. (See \cite{MR0442061} or \cite{MR603953} \S 1.2)

For each discrete valuation $v$ of the field $K=k(t)$ of rational functions over $k$, let $\pi_v$ be a uniformizing parameter and
$k_v = R_v/(\pi_v)$ be the residue field of the valuation ring $R_v=\{ r \in K | v(r) \ge 0 \}$.
Then we define the tame symbol $\partial_v : K^M_{l+1}(K) \rightarrow K^M_l(k_v)$ to be the epimorphism such that
$\partial ( \{ u_1, \dots, u_l, y \}) = v(y) \{\overline{u_1}, \dots, \overline{u_l} \}$ whenever $u_1, \dots, u_l$ are
units of the valuation ring $R_v$.

Let $v_\infty$ be the valuation on $K=k(t)$, which vanishes on $k$, such that $v_\infty (t) = -1$.
Every simple algebraic extension $L$ of $k$ is isomorphic to $k_v$ for some discrete valuation $v \ne v_\infty$ which corresponds to a
prime ideal $\mathfrak{p}$ of $k[t]$. The norm maps $N_v : K^M_l(k_v) \rightarrow K^M_l(k)$ are the unique homomorphisms such that,
for every $w \in K^M_{l+1}(k(t)),$ $\displaystyle \sum_{v} N_v \left( \partial_v w \right)=0$
where the sum is taken over all discrete valuations, including $v_\infty$ on $k(t)$, vanishing on $k$. This equality is called
the Weil reciprocity law. Note that we take $N_{v_\infty} = \rm{Id}$ for $v=v_\infty$.

Kato (\cite{MR603953} \S 1.7) has shown that these maps, if defined as compositions of norm maps for simple extensions for a given tower of
simple extensions, depend only on the field extension $L/k$, i.e., that it enjoys functoriality. See also \cite{MR689382}. It also enjoys
a projection formula similar to $(iii)$ of Lemma \ref{basic-norm}.
The following key lemma shows the compatibility between these two types of norm maps.

\begin{lemma} \label{norm-compatible}
For every finite field extension $k \subset L$, we have the following commutative diagram, where the vertical maps are the norm maps
and the horizontal maps are the homomorphisms in Proposition \ref{Milnor-map}:
$$\xymatrix{K^M_l (L) \ar[r]^-{\rho_l} \ar[d]^-{N_{L/k}} &
    H^{l}_{\M} \bigl(\Spec L , \Z(l) \bigr) \ar[d]^-{N_{L/k}} \\
    K^M_l (k) \ar[r]^-{\rho_l} &
    H^{l}_{\M} \bigl(\Spec k , \Z(l) \bigr)}$$
\end{lemma}
\begin{proof}
We will follow the same procedure which is used in \cite{MR2242284} for the proof.
Because of the functoriality property of the norm maps, we may assume that $[L:k]$ is a prime number $p$.
First, let us assume that $k$ has no extensions of degree prime to $p$. By Lemma (5.3) in \cite{MR0442061},
$K^M_l(L)$ is generated by the symbols of the form $x=\{x_1, \dots, x_{l-1}, y\}$ where $x_i \in k$ and $y \in L$. Then,
by the projection formula for Milnor's $K$-groups,
$\rho_l N_{L/k}\left(\{x_1, \dots, x_{l-1}, y\}\right) = \rho_l \left(\{x_1, \dots, x_{l-1}, N_{L/k}(y)\}\right)
            = \left( x_1, \dots, x_{l-1}, N_{L_k}(y) \right)$.
We also have $N_{L/k} \rho_l \left(\{x_1, \dots, x_{l-1}, y\}\right) = N_{L/k} \bigl( (x_1, \dots, x_{l-1}, y) \bigr)
=\left( x_1, \dots, x_{l-1}, N_{L_k}(y) \right)$ by $(iii)$ of Lemma \ref{basic-norm} and so we're done in this case.

Next, for the general case, let $k'$ be a maximal prime-to-$p$ extension of $k$. Then, by the previous case applied to $k'$ and by $(i)$ of Lemma
\ref{basic-norm}, we see that $z = N_{L/k} \rho_l(x) - \rho_l N_{L/k}(x)$, which is in the kernel of
$i_{k'/k}: H^{l}_{\M} \bigl(\Spec k , \Z(l) \bigr) \rightarrow H^{l}_{\M} \bigl(\Spec k' , \Z(l) \bigr)$, is
a torsion element of $H^{l}_{\M} \bigl(\Spec k , \Z(l) \bigr)$ of exponent prime to $p$.
In particular, if $L/k$ is a purely inseparable extension of degree $p$, then $y^p \in k$ and so $z$ is clearly killed by $p$, i.e., $z=0$.
Hence we may assume that $L/k$ is separable.

Since the kernel of $i_{L/k}: H^{l}_{\M} \bigl(\Spec k , \Z(l) \bigr) \rightarrow H^{l}_{\M} \bigl(\Spec L , \Z(l) \bigr)$ has exponent $p$,
it suffices to prove that $i_{L/k}(z) = 0$ to conclude $z=0$.
Now $L \otimes_k L$ is a finite product of fields $L_i$ with $[L_i:L] < p$ and we have the following commutative diagrams.
$$\xymatrix{K^M_l (L) \ar[r]^-{\oplus i_{L_i/L}} \ar[d]^-{N_{L/k}} &
    \oplus_i K^M_l (L_i) \ar[d]^-{\sum_i N_{L_i/L}} \\
    K^M_l (k) \ar[r]^{i_{L/k}} &
    K^M_l (L)}
\quad \quad  \xymatrix{H^{l}_{\M} \bigl(\Spec L , \Z(l) \bigr) \ar[r]^-{\oplus i_{L_i/L}} \ar[d]^-{N_{L/k}} &
    \oplus_i H^{l}_{\M} \bigl(\Spec L_i , \Z(l) \bigr) \ar[d]^-{\sum_i N_{L_i/L}} \\
    H^{l}_{\M} \bigl(\Spec k , \Z(l) \bigr) \ar[r]^-{i_{L/k}} &
    H^{l}_{\M} \bigl(\Spec L , \Z(l) \bigr)}
$$
The left diagram is the diagram (15) in p.387 of \cite{MR0442061} and the right diagram follows easily from
Definition \ref{def-norm}.
By induction on $p$, we have $i_{L/k}(z) =\oplus N_{L_i/L}\rho_l (i_{L_i/L}(x))- \oplus \rho_l N_{L_i/L} (i_{L_i/L}(x))=0$ and the proof is complete.
\end{proof}

\begin{lemma} \label{phi-map}
For any field $k$, there is a homomorphism $\phi_l : H^l_{\M} \bigl( \Spec k, \, \Z(l) \bigr) \rightarrow K^M_l(k)$ such that,
for each element $z \in H^{l}_{\M} \bigl(\Spec k , \Z(l) \bigr)$, there is an expression
$z = \displaystyle \sum_{j=1}^r N_{L_i/k} \left((\alpha_{1j}, \dots, \alpha_{lj})\right)$
where $L_1, \dots, L_r$ are finite field extensions of $k$,
$\alpha_{ij} \in GL_1(L_j) = L_j^\times$ ($1 \le i \le l$, $1 \le j \le r$) and an equality
$\displaystyle \phi_l (z) = \sum_j N_{L_j/k}\left( \{ \alpha_{1j}, \dots, \alpha_{lj} \} \right)$ in $K^M_l(k)$
\end{lemma}

\begin{proof}
For a tuple $z=(\theta_1, \theta_2, \dots, \theta_l)$, where $\theta_1, \theta_2, \dots, \theta_l$ are commuting matrices
in $GL_n(k)$, consider the vector space $E=k^n$ as an $R=k[t_1, t_1^{-1} \dots, t_l, t_l^{-1}]$-module,
on which $t_i$ acts as $\theta_i$. Since $E$ is of finite rank over $k$,
it has a composition series $0=E_0 \subset E_1 \subset \dots \subset E_r=E$
with simple factors $L_j = E_j / E_{j-1}$ ($j=1,\dots,r$).

Then, there exists a maximal ideal $\mathfrak{m}_j$ of $R$ such that $L_j \simeq R / \mathfrak{m}_j$.
So we see that $L_j$ is a finite field extension of $k$, and $\displaystyle z = \sum_{j=1}^r (\theta_1 | L_j, \dots, \theta_l | L_j)$,
where $\theta_i |L_j$ is the automorphism on $L_j$ induced by $\theta_i$.

Let us denote by $\alpha_{ij}$ the element of $L_j^\times$ which corresponds to $t_i$ (mod $\mathfrak{m}_j$) for $i=1, \dots, l$, then
$\displaystyle (\theta_1 | L_j, \dots, \theta_l | L_j) = N_{L_j/k} \left((\alpha_{1j}, \dots, \alpha_{lj})\right)$.

Since these factors $L_j$ are unique up to an order and a Milnor symbol vanishes if one of its coordinates is $1$,
the assignment $(\theta_1, \theta_2, \dots, \theta_l) \mapsto \sum_j N_{L_j/k} \left(\{\alpha_{1j}, \dots, \alpha_{lj}\}\right)$
gives us a well-defined homomorphism from $K_0 (k, \, \G_m^{\wedge l})$ to $K^M_l(k)$.

It remains to show that this homomorphism vanishes on $\partial K_0 (k \Delta^1, \, \G_m^{\wedge l})$.
Let $A_1(t), \dots, A_l(t)$ be commuting matrices in $GL_n(k[t])$, where $t$ is an indeterminate.
Then $M=k(t)^n$ can be considered as an $S=k(t)[t_1,t_1^{-1}, \dots, t_l, t_l^{-1}]$-module, on which $t_i$ acts as $A_i(t)$. Then
find a composition series $0=M_0 \subset M_1 \subset \dots \subset M_s=M$ with simple $S$-modules $Q_j = M_j / M_{j-1}$ ($j=1,\dots,s$)
and maximal ideals $\mathfrak{n}_j$ of $S$ such that $Q_j \simeq S / \mathfrak{n}_j$.
We also denote by $\beta_{ij}$ the element of $Q_j^\times$ which corresponds to $t_i$ (mod $\mathfrak{n}_j$) for $i=1, \dots, l$ and $j=1, \dots, s$.
Each $Q_j$ is a finite extension field of $k(t)$ and let $x = \sum_{j=1}^s N_{Q_j / k(t)} (\{\beta_{1j}, \dots, \beta_{lj}\}) \in K^M_l(k(t))$.
Now consider the element $\displaystyle y=\{x, {(t-1)/ t} \}$ in $K^M_{l+1}(k(t))$,
where the symbol $\{x, {(t-1)/ t} \}$ denotes $\sum_u \{x_{1u}, \dots, x_{lu}, {(t-1)/t} \}$ if $x = \sum_u \{x_{1u}, \dots, x_{lu}\}$ in $K^M_l(k(t))$.
Then $\partial_v (y) = - \phi_l \big( (A_1(0), \dots, A_l(0)) \big)$ if $\pi_v = t$ and
$\partial_v (y) = \phi_l \big( (A_1(1), \dots, A_l(1)) \big)$ if $\pi_v = t-1$.
Also, the image $\partial_v (y)$ is zero unless $v$ is the valuation associated with either $\pi_v = t-1$ or $\pi_v = t$.

Hence we have $\phi_l \big( (A_1(0), \dots, A_l(0)) \big) = \phi_l \big( (A_1(1), \dots, A_l(1)) \big)$ by the Weil reciprocity law for
the Milnor's $K$-groups.
\end{proof}

The isomorphism in the following theorem was first given by Nesterenko and Suslin (\cite{MR992981}) for Bloch's higher Chow groups.
Totaro, in \cite{MR1187705}, gave another proof of the theorem.
Suslin and Voevodsky, in Chapter 3 of \cite{MR1744945}, gave a proof of it for their motivic cohomology.
Here, we present another version of it for the Goodwillie-Lichtenbaum motivic complex
such that the isomorphism is given explicitly in the form which transforms
the multilinearity of the the symbols of Milnor into the corresponding properties of the symbols of Goodwillie and Lichtenbaum.

\begin{theorem} \label{Milnor-iso}
For any field $k$ and $l \ge 1$, the assignment $\{ a_1,a_2, \dots, a_l \} \mapsto (a_1,a_2,\dots,a_l) $ for each Steinberg symbol
$\{a_1,a_2, \dots, a_l \}$ gives rise to an isomorphism $K^M_l(k) \simeq H^l_{\M} \bigl( \Spec k, \, \Z(l) \bigr)$.
\end{theorem}
\begin{proof}
The case $l=1$ is straightforward and we assume $l \ge 2$. By Proposition \ref{Milnor-map}, the assignment $\{ a_1,a_2, \dots, a_l \} \mapsto (a_1,a_2,\dots,a_l) $
gives rise to a homomorphism $\rho_l$ from the Milnor's $K$-group $K^M_l(k)$ to the motivic cohomology group $H^l_{\M} \bigl( \Spec k, \, \Z(l) \bigr)$.

We also have a well-defined map $ \phi_l : H^l_{\M} \bigl( \Spec k, \, \Z(l) \bigr) \rightarrow K^M_l(k)$ in Lemma \ref{phi-map} and it suffices
to show that they are the inverses to each other.

It is clear that $\phi_l \circ \rho_l$ is the identity map on $K^M_l(k)$
since each Steinberg symbol is fixed by it. On the other hand, for each $z \in H^l_{\M} \bigl( \Spec k, \, \Z(l) \bigr)$,
$\displaystyle z = \sum_{j=1}^r N_{L_j/k} \left((\alpha_{1j}, \dots, \alpha_{lj})\right)$ for some finite field extensions
$L_1, \dots, L_r$ of $k$ and $\alpha_{ij} \in L_j$ ($1 \le i \le l$, $1 \le j \le r$). Then
$\displaystyle (\rho_l \circ \phi_l) (z) = \rho_l \left( \sum_j N_{L_j/k}\left( \{ \alpha_{1j}, \dots, \alpha_{lj} \} \right) \right)
                                =  \sum_j N_{L_j/k} \left( \rho_l \left( \{ \alpha_{1j}, \dots, \alpha_{lj} \} \right) \right)
                                =  \sum_j N_{L_j/k} \left( \left( \alpha_{1j}, \dots, \alpha_{lj} \right) \right) = z$ by Lemma \ref{norm-compatible}.
Therefore, $\rho_l \circ \phi_l$ is also the identity map and the proof is complete.
\end{proof}

\section{Multilinearity and Skew-symmetry for $H^{l-1}_{\mathcal{M}} \bigl( \rm{Spec} \, k, \, \Z(l) \bigr)$}

In \cite{MR2189214}, the author constructed a dilogarithm map $D: H^{1}_{\M} \bigl( \Spec k, \, \Z(2) \bigr) \rightarrow \R$ whenever $k$ is a subfield
of $\C$ such that $D$ satisfies certain bilinearity and skew-symmetry. (See Lemma 4.8 in \cite{MR2189214}).
Since $D$ can detect all the torsion-free elements of the motivic cohomology group , e.g.,
if $k$ is a number field (\cite{MR58:22016}, \cite{MR2001i:11082}),
we have expected that bilinearity and skew-symmetry for symbols should hold
for $D: H^{1}_{\M} \bigl( \Spec k, \, \Z(2) \bigr) \rightarrow \R$ in such cases.

In this section, we extend multilinearity and skew-symmetry results of the previous section to
the symbols in the motivic cohomology groups $H^{l-1}_{\M} \bigl( \Spec k, \, \Z(l) \bigr)$ when $k$ is a field.

$K_0 (k \Delta^{1}, \, \G_m^{\wedge l})$ ($l \ge 1$) can be identified with the abelian group
generated by $l$-tuples $(\theta_1,\dots,\theta_l)$ $\left( = \left( \theta_1(t),\dots,\theta_l(t)\right) \right)$
and certain explicit relations, where $\theta_1,\dots,\theta_l$ are commuting matrices in $GL_n(k [t])$ for various $n \ge 1$.
$K_0 (k \Delta^2, \, \G_m^{\wedge l})$ is identified with the abelian group
generated by the symbols $\left( \theta_1(x,y),\dots,\theta_l(x,y) \right)$ with commuting $\theta_1(x,y),\dots,\theta_l(x,y)
\in GL_n(k[x,y])$ and certain relations, and the boundary map $\partial$ on the motivic complex sends
$\left( \theta_1(x,y),\dots,\theta_l(x,y) \right)$ to
$\left( \theta_1(1-t,t),\dots,\theta_l(1-t,t) \right) - \left( \theta_1(0,t),\dots,\theta_l(0,t) \right)
+ \left( \theta_1(t,0),\dots,\theta_l(t,0) \right)$ in $K_0 (k \Delta^1, \, \G_m^{\wedge l})$.
The same symbol $(\theta_1,\dots,\theta_l)$ will denote the element in
$K_0 (k \Delta^1, \, \G_m^{\wedge l}) / \partial K_0 (k \Delta^2, \, \G_m^{\wedge l})$ represented by $(\theta_1,\dots,\theta_l)$,
by abuse of notation. The motivic cohomology group $H^{l-1}_{\M} \bigl( \Spec k, \, \Z(l) \bigr)$ is a subgroup of this quotient group,
which consists of the elements killed by $\partial$.

\begin{lemma} \label{simple1}
In $H^{l-1}_{\M} \bigl( \Spec k, \, \Z(l) \bigr)$, we have the following two simple relations of symbols for any commuting matrices
$\theta_1,\dots,\theta_l$ and any other commuting matrices $\psi_1,\dots,\psi_l$ in $GL_n(k[t])$:
\begin{align*}
-\left( \theta_1(t),\dots,\theta_l(t)\right) = \left( \theta_1(1-t),\dots,\theta_l(1-t)\right)
\end{align*}
$$\left( \theta_1(t),\dots,\theta_l(t) \right)+ \left( \psi_1(t),\dots,\psi_l(t) \right)
= \left( \theta_1(t)\oplus \psi_1(t), \dots, \theta_l(t) \oplus \psi_l(t) \right).$$
\end{lemma}

\begin{proof}
The second relation is immediate from definition of the motivic complex.
The first relation can be shown by applying the boundary map $\partial$ to the element $\left( \theta_1(x),\dots,\theta_l(x) \right)$
regarded as in $K_0 (k \Delta^2, \, \G_m^{\wedge l})$ and by noting that $(\theta_1,\dots,\theta_l) = 0$
in $H^{l-1}_{\M} \left( \Spec k, \, \Z(l) \right)$ when $\theta_1,\dots,\theta_l$ are constant matrices.
The fact that $(\theta_1,\dots,\theta_l) = 0$ for constant matrices $\theta_1,\dots,\theta_l$ is obtained simply by
applying the boundary map $\partial$ to the element $\left( \theta_1,\dots,\theta_l \right)$ regarded as in $K_0 (k \Delta^2, \, \G_m^{\wedge l})$.
\end{proof}

\begin{corollary}
Any element of the cohomology group $H^{l-1}_{\M} \bigl( \Spec k, \, \Z(l) \bigr)$ can be represented by a single expression $(\theta_1,\dots,\theta_l)$,
where $\theta_1,\dots,\theta_l$ are commuting matrices in $GL_n(k[t])$ for some nonnegative integer $n$.
\end{corollary}

We remark that the symbol $\left( \theta_1(t),\dots,\theta_l(t)\right)$ represents an element in $H^{l-1}_{\M} \bigl( \Spec k, \, \Z(l) \bigr)$
only when its image under the boundary map $\partial$ vanishes in $K_0 (k \Delta^0, \, \G_m^{\wedge l})$.

A tuple $\left( \theta_1(t),\dots,\theta_l(t)\right)$ where $\theta_1,\dots,\theta_l$ are commuting matrices in $GL_n(k[t])$
is called irreducible if $k[t]^n$ has no nontrivial proper submodule when regarded as a $k[t,x_1, x_1^{-1} \dots, x_l, x_l^{-1}]$-module
where $x_i$ acts on $k[t]^n$ via $\theta_i(t)$ for each $i=1, \dots, l$. Note that if $k(t)^n$ is regarded as
a $k(t)[x_1, x_1^{-1} \dots, x_l, x_l^{-1}]$-module with the same actions and if $M$ is a nontrivial proper
submodule of $k(t)^n$, then $M \cap k[t]^n$ is a nontrivial proper $k[t,x_1, x_1^{-1} \dots, x_l, x_l^{-1}]$-submodule of $k[t]^n$. Therefore,
$k(t)^n$ is irreducible as a $k(t)[x_1, x_1^{-1} \dots, x_l, x_l^{-1}]$-module if $\left( \theta_1(t),\dots,\theta_l(t)\right)$ is irreducible.

 It can be easily checked that, if two matrices $A, B \in GL_n(k)$ commute and $A$ is a block matrix of the form
$\displaystyle A = \begin{pmatrix} I & 0 \\ 0 & C \end{pmatrix}$ where $I$ is a matrix whose characteristic polynomial is a power of $x-1$ and
$C$ does not have 1 as an eigenvalue, then $B$ must be a block matrix $\displaystyle B = \begin{pmatrix} B_1 & 0 \\ 0 & B_2 \end{pmatrix}$,
where the blocks $B_1$ and $B_2$ are of compatible sizes with the blocks $I$ and $C$ of $A$. Therefore, we may easily relax the notion of
irreduciblity of a symbol $\left( \theta_1(t),\dots,\theta_l(t)\right)$ as an element in $K_0 (k \Delta^1, \, \G_m^{\wedge l})$
by declaring it \underline{irreducible} when its restriction to the largest submodule $V \subset k[t]^n$, where none of the restrictions of $\theta_1(t),\dots,\theta_l(t)$
has 1 as an eigenvalue, is irreducible.

\begin{theorem} \label{multilinear} (Multilinearity)
Suppose that $\varphi(t), \psi(t)$ and $\theta_1(t), \dots, \theta_l(t)$ (with $\theta_i(t)$ omitted) are commuting matrices
in $GL_n(k [t])$ such that the symbol represented by one of these matrices is irreducible in $K_0 (k \Delta^1, \, \G_m^{\wedge 1})$.
Assume further that the symbols $\bigl( \theta_1(t) ,\dots, \varphi(t),\dots, \theta_l(t) \bigr)$ and
$\bigl( \theta_1(t) ,\dots, \psi(t),\dots, \theta_l(t) \bigr)$ represent elements in $H^{l-1}_{\M} \bigl( \Spec k, \, \Z(l) \bigr)$.
Then $\bigl( \theta_1(t) ,\dots, \varphi(t) \psi(t),\dots, \theta_l(t) \bigr)$ represents an element in $H^{l-1}_{\M} \bigl( \Spec k, \, \Z(l) \bigr)$ and
$$\bigl( \theta_1(t) ,\dots, \varphi(t),\dots, \theta_l(t) \bigr) + \bigl( \theta_1(t) ,\dots, \psi(t),\dots, \theta_l(t) \bigr)
= \bigl( \theta_1(t) ,\dots, \varphi(t) \psi(t),\dots, \theta_l(t) \bigr)$$
in $H^{l-1}_{\M} \bigl( \Spec k, \, \Z(l) \bigr)$.
\end{theorem}

\begin{proof}
For simplicity of notation, we may assume that $i=1$ and prove the multilinearity on the first variable, i.e., we will want to show that
$$\bigl( \varphi(t),\theta_2(t), \dots, \theta_l(t) \bigr) + \bigl( \psi(t),\theta_2(t), \dots, \theta_l(t) \bigr)
= \bigl( \varphi(t)\psi(t),\theta_2(t), \dots, \theta_l(t) \bigr).$$
In this proof, all equalities are in $K_0 (k \Delta^{1}, \, \G_m^{\wedge l}) / \partial K_0 (k \Delta^{2}, \, \G_m^{\wedge l})$
unless mentioned otherwise and $1$ denotes the identity matrix $1_n$ of rank $n$ whenever appropriate.

Let $p(t)$ and $q(t)$ be matrices with entries in $k [t]$ such that $p(t)$ is invertible and
$p(t)$, $q(t)$ and $\theta_2(t), \dots, \theta_l(t)$ commute. Then
the boundary of the element
$$ \left( {\begin{pmatrix} 0 & 1 \\- p(y) & x y q(y)
\end{pmatrix}}, \ {\begin{pmatrix} \theta_2(y) & 0 \\ 0 & \theta_2(y) \end{pmatrix}}, \dots, {\begin{pmatrix} \theta_l(y) & 0 \\ 0 & \theta_l(y) \end{pmatrix}} \right)$$
of $K_0 (k \Delta^2, \, \G_m^{\wedge l})$ vanishes in
$H^{l-1}_{\M} \bigl( \Spec k, \, \Z(l) \bigr)$ by the definition of the cohomology group. Hence we have
\begin{multline*}
0 = \left( {\begin{pmatrix} 0 & 1 \\ -p(t) & (1-t) t q(t) \end{pmatrix}}, \
{\begin{pmatrix} \theta_2(t) & 0 \\ 0 & \theta_2(t) \end{pmatrix}}, \dots, {\begin{pmatrix} \theta_l(t) & 0 \\ 0 & \theta_l(t) \end{pmatrix}} \right) \\
- \left( {\begin{pmatrix} 0 & 1 \\ -p(t) & 0 \end{pmatrix}}, \
{\begin{pmatrix} \theta_2(t) & 0 \\ 0 & \theta_2(t) \end{pmatrix}}, \dots, {\begin{pmatrix} \theta_l(t) & 0 \\ 0 & \theta_l(t) \end{pmatrix}} \right)
+ \left( {\begin{pmatrix} 0 & 1 \\ p(0) & 0
\end{pmatrix}}, \  {\begin{pmatrix} \theta_2(0) & 0 \\ 0 & \theta_2(0) \end{pmatrix}}, \dots, {\begin{pmatrix} \theta_l(0) & 0 \\ 0 & \theta_l(0) \end{pmatrix}} \right).
\end{multline*}
But, as in the proof of Lemma \ref{simple1}, the last term, which is a tuple of constant matrices, is 0 and we have
\begin{multline} \label{bilinear-equ1}
\quad \quad \left( {\begin{pmatrix} 0 & 1 \\ -p(t) & (1-t) t q(t) \end{pmatrix}},\
{\begin{pmatrix} \theta_2(t) & 0 \\ 0 & \theta_2(t) \end{pmatrix}}, \dots, {\begin{pmatrix} \theta_l(t) & 0 \\ 0 & \theta_l(t) \end{pmatrix}} \right) \\
= \left( {\begin{pmatrix} 0 & 1 \\ -p(t) & 0 \end{pmatrix}}, \
{\begin{pmatrix} \theta_2(t) & 0 \\ 0 & \theta_2(t) \end{pmatrix}}, \dots, {\begin{pmatrix} \theta_l(t) & 0 \\ 0 & \theta_l(t) \end{pmatrix}} \right). \quad \quad
\end{multline}

Next, by taking the boundary of $\left( {\begin{pmatrix} 0 & 1 \\ -p(y) & (x+y) q(y)
\end{pmatrix}}, \ {\begin{pmatrix} \theta_2(y) & 0 \\ 0 & \theta_2(y) \end{pmatrix}}, \dots, {\begin{pmatrix} \theta_l(y) & 0 \\ 0 & \theta_l(y) \end{pmatrix}} \right)$, we get
\begin{multline} \label{bilinear-equ2}
\left( {\begin{pmatrix} 0 & 1 \\ - p(t) & q(t) \end{pmatrix}}, \
{\begin{pmatrix} \theta_2(t) & 0 \\ 0 & \theta_2(t) \end{pmatrix}}, \dots, {\begin{pmatrix} \theta_l(t) & 0 \\ 0 & \theta_l(t) \end{pmatrix}} \right) \\
= \left( {\begin{pmatrix} 0 & 1 \\ -p(t) & t q(t) \end{pmatrix}}, \
{\begin{pmatrix} \theta_2(t) & 0 \\ 0 & \theta_2(t) \end{pmatrix}}, \dots, {\begin{pmatrix} \theta_l(t) & 0 \\ 0 & \theta_l(t) \end{pmatrix}} \right) \\
- \left( {\begin{pmatrix} 0 & 1 \\ -p(0) & t q(0) \bigr)
\end{pmatrix}},\ {\begin{pmatrix} \theta_2(0) & 0 \\ 0 & \theta_2(0) \end{pmatrix}}, \dots, {\begin{pmatrix} \theta_l(0) & 0 \\ 0 & \theta_l(0) \end{pmatrix}} \right).
\end{multline}

If $p(t)$, $q(t)$ and $\theta_2(t), \dots, \theta_l(t)$ are replaced by $p(1-t)$,$(1-t)q(1-t)$ and $\theta_2(1-t), \dots, \theta_l(1-t)$ respectively in
(\ref{bilinear-equ2}), then we obtain
\begin{multline} \label{bilinear-equ3}
\left( {\begin{pmatrix} 0 & 1 \\ - p(1-t) & (1-t)q(1-t) \end{pmatrix}}, \
{\begin{pmatrix} \theta_2(1-t) & 0 \\ 0 & \theta_2(1-t) \end{pmatrix}}, \dots, {\begin{pmatrix} \theta_l(1-t) & 0 \\ 0 & \theta_l(1-t) \end{pmatrix}} \right) \\
= \left( {\begin{pmatrix} 0 & 1 \\ -p(1-t) & t(1-t) q(1-t) \end{pmatrix}}, \
{\begin{pmatrix} \theta_2(1-t) & 0 \\ 0 & \theta_2(1-t) \end{pmatrix}}, \dots, {\begin{pmatrix} \theta_l(1-t) & 0 \\ 0 & \theta_l(1-t) \end{pmatrix}} \right) \\
- \left( {\begin{pmatrix} 0 & 1 \\ -p(1) & t q(1) \bigr)
\end{pmatrix}},\ {\begin{pmatrix} \theta_2(1) & 0 \\ 0 & \theta_2(1) \end{pmatrix}}, \dots, {\begin{pmatrix} \theta_l(1) & 0 \\ 0 & \theta_l(1) \end{pmatrix}} \right).
\end{multline}

If we apply Lemma \ref{simple1} to the first term, the right hand side of the equality (\ref{bilinear-equ2}) can be written as
\begin{multline*}
- \left( {\begin{pmatrix} 0 & 1 \\ -p(1-t) & (1-t) q(1-t) \end{pmatrix}},\
{\begin{pmatrix} \theta_2(1-t) & 0 \\ 0 & \theta_2(1-t) \end{pmatrix}}, \dots, {\begin{pmatrix} \theta_l(1-t) & 0 \\ 0 & \theta_l(1-t) \end{pmatrix}} \right) \\
- \left( {\begin{pmatrix} 0 & 1 \\ -p(0) & t q(0) \bigr)
\end{pmatrix}},\  {\begin{pmatrix} \theta_2(0) & 0 \\ 0 & \theta_2(0) \end{pmatrix}}, \dots, {\begin{pmatrix} \theta_l(0) & 0 \\ 0 & \theta_l(0) \end{pmatrix}} \right).
\end{multline*}
By applying (\ref{bilinear-equ3}) to the first term and by (\ref{bilinear-equ2}), we have
\begin{align*}
&\quad \left( {\begin{pmatrix} 0 & 1 \\ -p(t) & q(t) \end{pmatrix}},
                    \ {\begin{pmatrix} \theta_2(t) & 0 \\ 0 & \theta_2(t) \end{pmatrix}}, \dots, {\begin{pmatrix} \theta_l(t) & 0 \\ 0 & \theta_l(t) \end{pmatrix}} \right) \\
&= - \left( {\begin{pmatrix} 0 & 1 \\ -p(1-t) & t(1-t) q(1-t) \end{pmatrix}},
                 \ {\begin{pmatrix} \theta_2(1-t) & 0 \\ 0 & \theta_2(1-t) \end{pmatrix}}, \dots, {\begin{pmatrix} \theta_l(1-t) & 0 \\ 0 & \theta_l(1-t) \end{pmatrix}} \right)  \\
&\hskip 1 in + \left( {\begin{pmatrix} 0 & 1 \\ -p(1) & t q(1) \bigr)
\end{pmatrix}} , {\begin{pmatrix} \theta_2(1) & 0 \\ 0 & \theta_2(1) \end{pmatrix}}, \dots, {\begin{pmatrix} \theta_l(1) & 0 \\ 0 & \theta_l(1) \end{pmatrix}} \right)\\
\displaybreak[0]
&\hskip 2 in - \left( {\begin{pmatrix} 0 & 1 \\ -p(0) & t q(0) \bigr)
\end{pmatrix}}, \  {\begin{pmatrix} \theta_2(0) & 0 \\ 0 & \theta_2(0) \end{pmatrix}}, \dots, {\begin{pmatrix} \theta_l(0) & 0 \\ 0 & \theta_l(0) \end{pmatrix}} \right) \\
&= \left( {\begin{pmatrix} 0 & 1 \\ -p(t) & t(1-t) q(t) \end{pmatrix}},
                \   {\begin{pmatrix} \theta_2(t) & 0 \\ 0 & \theta_2(t) \end{pmatrix}}, \dots, {\begin{pmatrix} \theta_l(t) & 0 \\ 0 & \theta_l(t) \end{pmatrix}} \right) \\
&\hskip 1 in + \left( {\begin{pmatrix} 0 & 1 \\ -p(1) & t q(1) \bigr)
\end{pmatrix}}, \  {\begin{pmatrix} \theta_2(1) & 0 \\ 0 & \theta_2(1) \end{pmatrix}}, \dots, {\begin{pmatrix} \theta_l(1) & 0 \\ 0 & \theta_l(1) \end{pmatrix}} \right) \\
\displaybreak[0]
&\hskip 2 in - \left( {\begin{pmatrix} 0 & 1 \\ -p(0) & t q(0) \bigr) \end{pmatrix}}, \
{\begin{pmatrix} \theta_2(0) & 0 \\ 0 & \theta_2(0) \end{pmatrix}}, \dots, {\begin{pmatrix} \theta_l(0) & 0 \\ 0 & \theta_l(0) \end{pmatrix}} \right) \\
&= \left( {\begin{pmatrix} 0 & 1 \\ -p(t) & 0 \end{pmatrix}}, \
{\begin{pmatrix} \theta_2(t) & 0 \\ 0 & \theta_2(t) \end{pmatrix}}, \dots, {\begin{pmatrix} \theta_l(t) & 0 \\ 0 & \theta_l(t) \end{pmatrix}} \right) \\
&\hskip 1 in + \left( {\begin{pmatrix} 0 & 1 \\ -p(1) & t q(1)
\end{pmatrix}}, \ {\begin{pmatrix} \theta_2(1) & 0 \\ 0 & \theta_2(1) \end{pmatrix}}, \dots, {\begin{pmatrix} \theta_l(1) & 0 \\ 0 & \theta_l(1) \end{pmatrix}} \right) \\
&\hskip 2 in - \left( {\begin{pmatrix} 0 & 1 \\ -p(0) & t q(0) \bigr)\end{pmatrix}}, \
{\begin{pmatrix} \theta_2(0) & 0 \\ 0 & \theta_2(0) \end{pmatrix}}, \dots, {\begin{pmatrix} \theta_l(0) & 0 \\ 0 & \theta_l(0) \end{pmatrix}} \right).
\end{align*}

The second equality is obtained by applying Lemma \ref{simple1} to the first term and the last equality is by (\ref{bilinear-equ1}).
Now by setting $p(t) = \varphi(t)\psi(t)$ and $q(t) = \varphi(t) + \psi(t)$ in the above equality, we have
\begin{multline} \label{bilinear-equ4}
\left( {\begin{pmatrix} 0 & 1 \\ -\varphi(t)\psi(t) & \varphi(t)+ \psi(t) \end{pmatrix}},
                    \ {\begin{pmatrix} \theta_2(t) & 0 \\ 0 & \theta_2(t) \end{pmatrix}}, \dots, {\begin{pmatrix} \theta_l(t) & 0 \\ 0 & \theta_l(t) \end{pmatrix}} \right) \\
= \left( {\begin{pmatrix} 0 & 1 \\ -\varphi(t)\psi(t) & 0 \end{pmatrix}}, \ {\begin{pmatrix} \theta_2(t) & 0 \\ 0 & \theta_2(t) \end{pmatrix}}, \dots, {\begin{pmatrix} \theta_l(t) & 0 \\ 0 & \theta_l(t) \end{pmatrix}} \right) \hskip 1.5 in \\
+ \left( {\begin{pmatrix} 0 & 1 \\ -\varphi(1)\psi(1) & t \bigl( \varphi(1)+ \psi(1) \bigr)
\end{pmatrix}}, \ {\begin{pmatrix} \theta_2(1) & 0 \\ 0 & \theta_2(1) \end{pmatrix}}, \dots, {\begin{pmatrix} \theta_l(1) & 0 \\ 0 & \theta_l(1) \end{pmatrix}} \right) \\
- \left( {\begin{pmatrix} 0 & 1 \\ -\varphi(0)\psi(0) & t \bigl( \varphi(0)+ \psi(0) \bigr) \end{pmatrix}}, \
{\begin{pmatrix} \theta_2(0) & 0 \\ 0 & \theta_2(0) \end{pmatrix}}, \dots, {\begin{pmatrix} \theta_l(0) & 0 \\ 0 & \theta_l(0) \end{pmatrix}} \right).
\end{multline}

Similarly, with $p(t) = \varphi(t)\psi(t)$ and $q(t) = 1 + \varphi(t)\psi(t)$ this time, we get
\begin{multline} \label{bilinear-equ5}
\left( {\begin{pmatrix} 0 & 1 \\ -\varphi(t)\psi(t) & 1 + \varphi(t)\psi(t) \end{pmatrix}},
                    \ {\begin{pmatrix} \theta_2(t) & 0 \\ 0 & \theta_2(t) \end{pmatrix}}, \dots, {\begin{pmatrix} \theta_l(t) & 0 \\ 0 & \theta_l(t) \end{pmatrix}} \right) \\
= \left( {\begin{pmatrix} 0 & 1 \\ -\varphi(t)\psi(t) & 0 \end{pmatrix}}, \ {\begin{pmatrix} \theta_2(t) & 0 \\ 0 & \theta_2(t) \end{pmatrix}}, \dots, {\begin{pmatrix} \theta_l(t) & 0 \\ 0 & \theta_l(t) \end{pmatrix}} \right) \hskip 1.5 in \\
+ \left( {\begin{pmatrix} 0 & 1 \\ -\varphi(1)\psi(1) & t \bigl( 1 + \varphi(1)\psi(1) \bigr)
\end{pmatrix}}, \ {\begin{pmatrix} \theta_2(1) & 0 \\ 0 & \theta_2(1) \end{pmatrix}}, \dots, {\begin{pmatrix} \theta_l(1) & 0 \\ 0 & \theta_l(1) \end{pmatrix}} \right) \\
- \left( {\begin{pmatrix} 0 & 1 \\ -\varphi(0)\psi(0) & t \bigl( 1 + \varphi(0)\psi(0) \bigr) \end{pmatrix}}, \
{\begin{pmatrix} \theta_2(0) & 0 \\ 0 & \theta_2(0) \end{pmatrix}}, \dots, {\begin{pmatrix} \theta_l(0) & 0 \\ 0 & \theta_l(0) \end{pmatrix}} \right).
\end{multline}

The first terms on the right of (\ref{bilinear-equ4}) and (\ref{bilinear-equ5}) are the same,
so by subtracting (\ref{bilinear-equ5}) from (\ref{bilinear-equ4}), we obtain

\begin{align*}
&\left( {\begin{pmatrix} 0 & 1 \\ -\varphi(t)\psi(t) & \varphi(t)+ \psi(t) \end{pmatrix}},
                    \ {\begin{pmatrix} \theta_2(t) & 0 \\ 0 & \theta_2(t) \end{pmatrix}}, \dots, {\begin{pmatrix} \theta_l(t) & 0 \\ 0 & \theta_l(t) \end{pmatrix}} \right)\\
&\hskip 1 in - \left( {\begin{pmatrix} 0 & 1 \\ -\varphi(t)\psi(t) & 1 + \varphi(t)\psi(t) \end{pmatrix}},
                    \ {\begin{pmatrix} \theta_2(t) & 0 \\ 0 & \theta_2(t) \end{pmatrix}}, \dots, {\begin{pmatrix} \theta_l(t) & 0 \\ 0 & \theta_l(t) \end{pmatrix}} \right) \\
&=\left( {\begin{pmatrix} 0 & 1 \\ -\varphi(1)\psi(1) & t \bigl( \varphi(1)+ \psi(1) \bigr)
\end{pmatrix}}, \ {\begin{pmatrix} \theta_2(1) & 0 \\ 0 & \theta_2(1) \end{pmatrix}}, \dots, {\begin{pmatrix} \theta_l(1) & 0 \\ 0 & \theta_l(1) \end{pmatrix}} \right) \\
&\hskip 1 in - \left( {\begin{pmatrix} 0 & 1 \\ -\varphi(0)\psi(0) & t \bigl( \varphi(0)+ \psi(0) \bigr) \end{pmatrix}}, \
{\begin{pmatrix} \theta_2(0) & 0 \\ 0 & \theta_2(0) \end{pmatrix}}, \dots, {\begin{pmatrix} \theta_l(0) & 0 \\ 0 & \theta_l(0) \end{pmatrix}} \right) \\
&\hskip 0.5 in - \left( {\begin{pmatrix} 0 & 1 \\ -\varphi(1)\psi(1) & t \bigl( 1 + \varphi(1)\psi(1) \bigr)
\end{pmatrix}}, \ {\begin{pmatrix} \theta_2(1) & 0 \\ 0 & \theta_2(1) \end{pmatrix}}, \dots, {\begin{pmatrix} \theta_l(1) & 0 \\ 0 & \theta_l(1) \end{pmatrix}} \right) \\
&\hskip 1.5 in + \left( {\begin{pmatrix} 0 & 1 \\ -\varphi(0)\psi(0) & t \bigl( 1 + \varphi(0)\psi(0) \bigr) \end{pmatrix}}, \
{\begin{pmatrix} \theta_2(0) & 0 \\ 0 & \theta_2(0) \end{pmatrix}}, \dots, {\begin{pmatrix} \theta_l(0) & 0 \\ 0 & \theta_l(0) \end{pmatrix}} \right).
\end{align*}

Now we state our claim:
\vskip 0.1 in
\bf{Claim: } \it{The right hand side of the above equality is equal to 0.}
\vskip 0.1 in
\rm Once the claim is proved, we obtain the following equality.
\begin{multline} \label{bilinear-equ6}
\quad \left( {\begin{pmatrix} 0 & 1 \\ -\varphi(t)\psi(t) & \varphi(t)+ \psi(t) \end{pmatrix}},
                    \ {\begin{pmatrix} \theta_2(t) & 0 \\ 0 & \theta_2(t) \end{pmatrix}}, \dots, {\begin{pmatrix} \theta_l(t) & 0 \\ 0 & \theta_l(t) \end{pmatrix}} \right)\\
= \left( {\begin{pmatrix} 0 & 1 \\ -\varphi(t)\psi(t) & 1 + \varphi(t)\psi(t) \end{pmatrix}},
                    \ {\begin{pmatrix} \theta_2(t) & 0 \\ 0 & \theta_2(t) \end{pmatrix}}, \dots, {\begin{pmatrix} \theta_l(t) & 0 \\ 0 & \theta_l(t) \end{pmatrix}} \right). \quad
\end{multline}
To prove the claim, note first that, by our assumption, one of  $\varphi(t)$, $\psi(t)$ and $\theta_2(t), \dots, \theta_l(t)$, denoted $\theta(t)$,
is irreducible on the largest submodule $V \subset k[t]^n$, where none of the restrictions of $\theta_1(t),\dots,\theta_l(t)$
has 1 as an eigenvalue. We may easily assume that $V=k[t]^n$ since all the symbols under our interest vanish on the complement of $V$ in $k[t]^n$.
Then all of $\varphi(t)$, $\psi(t)$ and $\theta_2(t), \dots, \theta_l(t)$ can be written as polynomials of $\theta(t)$ with coefficients in $k(t)$.
Since $\bigl( \varphi(0), \theta_2(0),\dots,\theta_l(0)\bigr) = \bigl( \varphi(1), \theta_2(1),\dots,\theta_l(1)\bigr)$
in $K_0 (k, \, \G_m^{\wedge l})$ by our assumption, it follows that $S \varphi(0) S^{-1} = \varphi(1)$, $S \theta_i(0) S^{-1}=\theta_i(1)$
for every legitimate $i$, for some $S \in GL_n(k)$. Now, it is immediate that, in $K_0 (k \Delta^1, \, \G_m^{\wedge l})$,
\begin{align*}
&\left( {\begin{pmatrix} 0 & 1 \\ -\varphi(1)\psi(1) & t \bigl( \varphi(1)+ \psi(1) \bigr)
\end{pmatrix}}, \ {\begin{pmatrix} \theta_2(1) & 0 \\ 0 & \theta_2(1) \end{pmatrix}}, \dots, {\begin{pmatrix} \theta_l(1) & 0 \\ 0 & \theta_l(1) \end{pmatrix}} \right) \\
&\hskip 1 in = \left( {\begin{pmatrix} 0 & 1 \\ -\varphi(0)\psi(0) & t \bigl( \varphi(0)+ \psi(0) \bigr) \end{pmatrix}}, \
{\begin{pmatrix} \theta_2(0) & 0 \\ 0 & \theta_2(0) \end{pmatrix}}, \dots, {\begin{pmatrix} \theta_l(0) & 0 \\ 0 & \theta_l(0) \end{pmatrix}} \right) \\
&\text{and }  \left( {\begin{pmatrix} 0 & 1 \\ -\varphi(1)\psi(1) & t \bigl( 1 + \varphi(1)\psi(1) \bigr)
\end{pmatrix}}, \ {\begin{pmatrix} \theta_2(1) & 0 \\ 0 & \theta_2(1) \end{pmatrix}}, \dots, {\begin{pmatrix} \theta_l(1) & 0 \\ 0 & \theta_l(1) \end{pmatrix}} \right) \\
&\hskip 1 in = \left( {\begin{pmatrix} 0 & 1 \\ -\varphi(0)\psi(0) & t \bigl( 1 + \varphi(0)\psi(0) \bigr) \end{pmatrix}}, \
{\begin{pmatrix} \theta_2(0) & 0 \\ 0 & \theta_2(0) \end{pmatrix}}, \dots, {\begin{pmatrix} \theta_l(0) & 0 \\ 0 & \theta_l(0) \end{pmatrix}} \right).
\end{align*}
Therefore, the proof of the claim is complete.

Thanks to the identities (\ref{linear-equa}) and (\ref{linear-equb}), we have, by (\ref{bilinear-equ6}),
\begin{align*}
&\quad \ \left( {\begin{pmatrix} \psi(t) & 1 \\ 0 & \varphi(t) \end{pmatrix}}, \ {\begin{pmatrix} \theta_2(t) & 0 \\ 0 & \theta_2(t) \end{pmatrix}}, \dots, {\begin{pmatrix} \theta_l(t) & 0 \\ 0 & \theta_l(t) \end{pmatrix}} \right) \\
\displaybreak[0]
&=\left( {\begin{pmatrix} 0 & 1 \\ -\varphi(t) \psi(t) & \varphi(t)+\psi(t) \end{pmatrix}}, \ {\begin{pmatrix} \theta_2(t) & 0 \\ 0 & \theta_2(t) \end{pmatrix}}, \dots, {\begin{pmatrix} \theta_l(t) & 0 \\ 0 & \theta_l(t) \end{pmatrix}} \right) \\
\displaybreak[0]
&=\left( {\begin{pmatrix} 0 & 1 \\ -\varphi(t)\psi(t) & 1+\varphi(t)\psi(t) \end{pmatrix}}, \ {\begin{pmatrix} \theta_2(t) & 0 \\ 0 & \theta_2(t) \end{pmatrix}}, \dots, {\begin{pmatrix} \theta_l(t) & 0 \\ 0 & \theta_l(t) \end{pmatrix}} \right) \\
&=\left( {\begin{pmatrix} 1 & 1 \\ 0 & \varphi(t)\psi(t) \end{pmatrix}}, \ {\begin{pmatrix} \theta_2(t) & 0 \\ 0 & \theta_2(t) \end{pmatrix}}, \dots, {\begin{pmatrix} \theta_l(t) & 0 \\ 0 & \theta_l(t) \end{pmatrix}} \right).
\end{align*}
Hence $\bigl( \varphi(t),\,\theta_2(t), \dots, \theta_l(t) \bigr) + \bigl( \psi(t), \,\theta_2(t), \dots, \theta_l(t) \bigr)
= \bigl( \varphi(t)\psi(t), \, \theta_2(t), \dots, \theta_l(t) \bigr)$, as required.
\end{proof}

The irreducibility assumption in Theorem \ref{multilinear} is used only to justify the claim in the proof of the theorem.
The various conditions in the following corollary can replace the irreducibility assumption in the theorem.
We state the multilinearity of symbols only in the
first coordinate to simplify the notation, but a similar statement in another coordinate holds obviously.

\begin{corollary} \label{multilinear1}
Suppose that $\varphi(t), \psi(t), \theta_2(t), \dots, \theta_l(t)$ are commuting matrices in $GL_n(k [t])$
and that the symbols $\bigl( \varphi(t), \theta_2(t) ,\dots, \theta_l(t) \bigr)$ and
$\bigl( \psi(t), \theta_2(t) ,\dots, \theta_l(t) \bigr)$ represent elements in $H^{l-1}_{\M} \bigl( \Spec k, \, \Z(l) \bigr)$.
Then $\bigl( \varphi(t)\psi(t), \theta_2(t) ,\dots, \theta_l(t) \bigr)$ represents an element in $H^{l-1}_{\M} \bigl( \Spec k, \, \Z(l) \bigr)$ and
$$\bigl( \varphi(t), \theta_2(t) ,\dots, \theta_l(t) \bigr) + \bigl( \psi(t), \theta_2(t) ,\dots, \theta_l(t) \bigr)
= \bigl( \varphi(t)\psi(t), \theta_2(t) ,\dots, \theta_l(t) \bigr)$$
in $H^{l-1}_{\M} \bigl( \Spec k, \, \Z(l) \bigr)$ if one of the following assumptions is satisfied:

(i) The symbol $\bigl( \varphi(t), \psi(t), \theta_2(t) ,\dots, \theta_l(t) \bigr)$ is irreducible and
$k$ is a field of characteristic $0$ or $n < char(k)$.

(ii) There exists a filtration $0=V_0 \subset V_1 \subset \dots \subset V_n=k[t]^n$ of $k[t,x_0, x_0^{-1} \dots, x_l, x_l^{-1}]$-modules where
$x_0$ and $x_1$ act via $\varphi(t)$ and $\psi(t)$ and $x_i$ acts via $\theta_i(t)$ for $i \ge 2$ such that the restriction of the symbol
$\bigl( \varphi(t), \psi(t), \theta_2(t) ,\dots, \theta_l(t) \bigr)$ to each $V_{i+1}/V_i$ ($i=0, \dots, n-1$) is irreducible and
$k$ is of characteristic $0$ or $n < char(k)$.

(iii) One of the matrices $\varphi(t), \psi(t), \theta_2(t) ,\dots, \theta_l(t)$ has a characteristic polynomial equal to its minimal polynomial.
This is the case, for example, when one of the matrices is a companion matrix of a polynomial with coefficients in $k[t]$ and constant term in $k^\times$.
\end{corollary}

\begin{proof}
$(i)$ $k(t)^n$ as a $k(t)[x_0, x_0^{-1} \dots, x_l, x_l^{-1}]$-module,
where $x_0$ and $x_1$ act via $\varphi(t)$ and $\psi(t)$ and $x_i$ acts via $\theta_i(t)$ for $i \ge 2$, is irreducible. Therefore, it is a field
extension of $k(t)$ of degree $n$. By our assumption on the field $k$, it is generated by a primitive element, say $\theta(t)$, and all of
$\varphi(t), \psi(t), \theta_2(t) ,\dots, \theta_l(t)$ can be written as polynomials of $\theta(t)$ with coefficients in $k(t)$.
So the claim in the proof of Theorem \ref{multilinear} holds and we obtain the multilinearity.

$(ii)$ is an obvious consequence of $(i)$.

$(iii)$ it true since any matrix which commutes with a given companion matrix of a polynomial can be written as a polynomial of the companion matrix.
(Theorem 5 of Chapter 1 in \cite{MR0201472})
\end{proof}

In the following corollary, we don't require the commutativity of $\varphi(t)$ and $\psi(t)$.

\begin{corollary} \label{multilinear2}
Suppose that $\theta_2(t), \dots, \theta_l(t)$ are commuting matrices in $GL_n(k [t])$
which commute also with $\varphi(t), \psi(t) \in GL_n(k[t])$ and that the symbols $\bigl( \varphi(t), \theta_2(t) ,\dots, \theta_l(t) \bigr)$ and
$\bigl( \psi(t), \theta_2(t) ,\dots, \theta_l(t) \bigr)$ represent elements in $H^{l-1}_{\M} \bigl( \Spec k, \, \Z(l) \bigr)$.
Then $\bigl( \varphi(t)\psi(t), \theta_2(t) ,\dots, \theta_l(t) \bigr)$ represents an element in $H^{l-1}_{\M} \bigl( \Spec k, \, \Z(l) \bigr)$ and
$$\bigl( \varphi(t), \theta_2(t) ,\dots, \theta_l(t) \bigr) + \bigl( \psi(t), \theta_2(t) ,\dots, \theta_l(t) \bigr)
= \bigl( \varphi(t)\psi(t), \theta_2(t) ,\dots, \theta_l(t) \bigr)$$
in $H^{l-1}_{\M} \bigl( \Spec k, \, \Z(l) \bigr)$ if one of the following assumptions is satisfied:

(i) $\varphi(0)=\varphi(1)$, $\psi(0)=\psi(1)$ and $\theta_i(0)=\theta_i(1)$ for $i=2, \dots, l$ as matrices in $GL_n(k)$.

(ii) $\theta_i(0)$ or $\theta_i(1)$ has $n$ distinct eigenvalues for some $i=2, \dots, l$.
\end{corollary}
\begin{proof}
$(i)$ clearly guarantees the claim in the proof of Theorem \ref{multilinear}.

$(ii)$ We may assume that none of $\varphi(0), \psi(0), \theta_2(0) ,\dots, \theta_l(0)$ has 1 as an eigenvalue.
If $\theta_i(0)$ has $n$ distinct eigenvalues for some $i$, then $\theta_i(1)$ also has the same $n$ distinct eigenvalues since
$(\theta_i(0)) = (\theta_i(1))$ in $K_0 (k, \, \G_m^{\wedge 1})$ by the assumption that $\bigl( \varphi(t), \theta_2(t) ,\dots, \theta_l(t) \bigr)$
belongs to $H^{l-1}_{\M} \bigl( \Spec k, \, \Z(l) \bigr)$.
Also, each of $\varphi(0), \psi(0), \theta_2(0) ,\dots, \theta_l(0)$ are
diagonalizable by the same similarity matrix by the commutativity of the matrices with $\theta_i(0)$.
Let us denote the tuples of joint eigenvalues by $(a_i, b_i, c_{2i}, \dots, c_{li})$ for $i=1, \dots, n$.
A similar statement is true for $\varphi(1), \psi(1), \theta_2(1) ,\dots, \theta_l(1)$ and their joint eigenvalues are denoted by
$(a'_i, b'_i, c'_{2i}, \dots, c'_{li})$ for $i=1, \dots, n$. By permuting the indices $i$ if necessary, we may assume that
$a_i = a'_i$, $b_i = b'_i$, $c_{ji}=c'_{ji}$ for $j=2, \dots, l$ and $i=1, \dots, n$.
Then the claim in the proof of Theorem \ref{multilinear} holds since
\begin{align*}
&\left( {\begin{pmatrix} 0 & 1 \\ -\varphi(1)\psi(1) & t \bigl( \varphi(1)+ \psi(1) \bigr)
\end{pmatrix}}, \ {\begin{pmatrix} \theta_2(1) & 0 \\ 0 & \theta_2(1) \end{pmatrix}}, \dots, {\begin{pmatrix} \theta_l(1) & 0 \\ 0 & \theta_l(1) \end{pmatrix}} \right) \\
&\hskip 0.3 in = \sum_{i=1}^n \left( {\begin{pmatrix} 0 & 1 \\ -a_i b_i & t \bigl( a_i+ b_i \bigr)
\end{pmatrix}}, \ {\begin{pmatrix} c_{2i} & 0 \\ 0 & c_{2i} \end{pmatrix}}, \dots, {\begin{pmatrix} c_{li} & 0 \\ 0 & c_{li} \end{pmatrix}} \right) \\
&\hskip 1 in = \left( {\begin{pmatrix} 0 & 1 \\ -\varphi(0)\psi(0) & t \bigl( \varphi(0)+ \psi(0) \bigr) \end{pmatrix}}, \
{\begin{pmatrix} \theta_2(0) & 0 \\ 0 & \theta_2(0) \end{pmatrix}}, \dots, {\begin{pmatrix} \theta_l(0) & 0 \\ 0 & \theta_l(0) \end{pmatrix}} \right) \\
&\text{and similarly}  \left( {\begin{pmatrix} 0 & 1 \\ -\varphi(1)\psi(1) & t \bigl( 1 + \varphi(1)\psi(1) \bigr)
\end{pmatrix}}, \ {\begin{pmatrix} \theta_2(1) & 0 \\ 0 & \theta_2(1) \end{pmatrix}}, \dots, {\begin{pmatrix} \theta_l(1) & 0 \\ 0 & \theta_l(1) \end{pmatrix}} \right) \\
&\hskip 1 in = \left( {\begin{pmatrix} 0 & 1 \\ -\varphi(0)\psi(0) & t \bigl( 1 + \varphi(0)\psi(0) \bigr) \end{pmatrix}}, \
{\begin{pmatrix} \theta_2(0) & 0 \\ 0 & \theta_2(0) \end{pmatrix}}, \dots, {\begin{pmatrix} \theta_l(0) & 0 \\ 0 & \theta_l(0) \end{pmatrix}} \right).
\end{align*}

\end{proof}

\begin{note}
(i) In Theorem \ref{multilinear}, the commutativity of $\varphi(t)$ and $\psi(t)$ would not have been necessary if we wanted merely to define the symbols
$\bigl( \varphi(t),\theta_2(t), \dots, \theta_l(t) \bigr)$ and $\bigl( \psi(t),\theta_2(t), \dots, \theta_l(t) \bigr)$.
But, if we do not insist the commutativity of these two matrices, then $\bigl( \varphi(t)\psi(t),\theta_2(t), \dots, \theta_l(t) \bigr)$
does not have to represent an element in $H^{l-1}_{\M} \bigl( \Spec k, \, \Z(l) \bigr)$ even if the symbols
$\bigl( \varphi(t),\theta_2(t), \dots, \theta_l(t) \bigr)$ and $\bigl( \psi(t),\theta_2(t), \dots, \theta_l(t) \bigr)$ do.

For example, take $l=2$ and let $a, b \in k-\{0,1\}$ be two distinct numbers and take any $c \in k-\{0,1\}$. Let
$$\varphi(t) = \begin{pmatrix} (a+b)t & {\frac {(a+b)^2} {ab}} t(1-t) -1 \\ ab & (a+b)(1-t) \end{pmatrix},
\psi(t) = \begin{pmatrix} a & 0 \\ 0 & b \end{pmatrix},
\theta(t) = \begin{pmatrix} c & 0 \\ 0 & c \end{pmatrix} $$
Then the boundaries of both $\bigl( \psi(t),\, \theta(t) \bigr)$ and $\bigl( \varphi(t),\, \theta(t) \bigr)$ are 0, but
the boundary of $\bigl(\varphi(t) \psi(t),\, \theta(t) \bigr)$ is not 0 in $K_0 (k \Delta^0, \, \G_m^{\wedge 2})$.

(ii) The irreducibility condition in Theorem \ref{multilinear} or other similar assumptions in Corollary \ref{multilinear1} and \ref{multilinear2}
are essential. For example, take $l=1$ and let $a, b \in k - \{0,1\}$ be two distinct elements.
Find any distinct $c, d \in k-\{0,\pm 1\}$ such that the set
$\{a, acd, bc, bd \}$ is not equal to $\{ac, ad, b, acd\}$. Consider
$$A(t) = \begin{pmatrix} a & 0& 0& 0 \\ 0 & a & 0 & 0 \\ 0 & 0 & b & 0 \\ 0 & 0 & 0 & b \end{pmatrix}, \
B(t) = \begin{pmatrix} 0 & -cd & 0& 0 \\ 1 & (c+d)t + (1+cd)(1-t) & 0 & 0 \\ 0 & 0 & 0 & -cd \\ 0 & 0 & 1 & (c+d)(1-t) + (1+cd)t \end{pmatrix}.$$
Then $A(0)=A(1)$ and $(B(0)) = (1) + (cd) + (c) + (d) = (B(1))$ in $K_0 (k, \, \G_m^{\wedge 1})$.
But, $(A(0)B(0)) = (a) + (acd) + (bc) + (bd) \ne (ac) + (ad) + (b) + (bcd) = (A(1)B(1))$ in $K_0 (k, \, \G_m^{\wedge 1})$ and thus
$(A(t)B(t))$ does not represent an element in $H^{0}_{\M} \bigl( \Spec k, \, \Z(1) \bigr)$
\end{note}

\begin{proposition} \label{skewsymmetry}
(Skew-Symmetry) Suppose that $\theta_1(t),\dots, \theta_l(t) \in GL_n(k [t])$ commute and one of the symbols represented by
$\theta_1(t),\dots, \theta_{l-1}(t)$ or $\theta_l(t)$ is irreducible.
If $\bigl( \theta_1(t),\dots, \theta_l(t) \bigr)$ represents an element in $H^{l-1}_{\M} \bigl( \Spec k, \, \Z(l) \bigr)$ ($l \ge 2$), then
$\bigl( \theta_1(t),\dots,\theta_i(t),\dots,\theta_j(t),\dots,\theta_l(t) \bigr) =
-\bigl( \theta_1(t),\dots,\theta_j(t),\dots,\theta_i(t),\dots,\theta_l(t) \bigr)$ in $H^{l-1}_{\M} \bigl( \Spec k, \, \Z(l) \bigr)$.
\end{proposition}

\begin{proof}
For simplicity of notations, we assume that $i=1$ and $j=2$. Let $\varphi=\theta_1$ and $\psi=\theta_2$.
An argument similar to the one utilized in the proof of Theorem \ref{multilinear} can be used to prove that
\begin{multline} \label{skewsymmetryeq}
\left( {\begin{pmatrix} 0 & 1 \\ -\varphi(t) \psi(t) & \varphi(t)+\psi(t) \end{pmatrix}},
        \ {\begin{pmatrix} 0 & 1 \\ -\varphi(t) \psi(t) & \varphi(t)+\psi(t) \end{pmatrix}},\
        {\begin{pmatrix} \theta_3(t) & 0 \\ 0 & \theta_3(t) \end{pmatrix}}, \dots, {\begin{pmatrix} \theta_l(t) & 0 \\ 0 & \theta_l(t) \end{pmatrix}} \right) \\
=\left( {\begin{pmatrix} 0 & 1 \\ -\varphi(t)\psi(t) & 1+\varphi(t)\psi(t) \end{pmatrix}},
        \ {\begin{pmatrix}  0 & 1 \\ -\varphi(t)\psi(t) & 1+\varphi(t)\psi(t) \end{pmatrix}},\
        {\begin{pmatrix} \theta_3(t) & 0 \\ 0 & \theta_3(t) \end{pmatrix}}, \dots, {\begin{pmatrix} \theta_l(t) & 0 \\ 0 & \theta_l(t) \end{pmatrix}} \right).
\end{multline}
\begin{multline*}
\text{ Just replace }\left( {\begin{pmatrix} 0 & 1 \\- p(t) &  q(t) \end{pmatrix}}, \
{\begin{pmatrix} \theta_2(t) & 0 \\ 0 & \theta_2(t) \end{pmatrix}}, \dots, {\begin{pmatrix} \theta_l(t) & 0 \\ 0 & \theta_l(t) \end{pmatrix}} \right) \\
\text{ by }\left( {\begin{pmatrix} 0 & 1 \\- p(t) & q(t) \end{pmatrix}}, \
{\begin{pmatrix} 0 & 1 \\- p(t) &  q(t) \end{pmatrix}}, \
{\begin{pmatrix} \theta_3(t) & 0 \\ 0 & \theta_3(t) \end{pmatrix}}, \dots, {\begin{pmatrix} \theta_l(t) & 0 \\ 0 & \theta_l(t) \end{pmatrix}} \right)
\end{multline*}
and make similar replacements throughout the course of the proof of the claim in the proof of Theorem \ref{multilinear}.
Then note that
\begin{align*}
&\left( {\begin{pmatrix} 0 & 1 \\ -\varphi(0)\psi(0) & t \bigl( \varphi(0)+ \psi(0) \bigr) \end{pmatrix}}, \
 {\begin{pmatrix} 0 & 1 \\ -\varphi(0)\psi(0) & t \bigl( \varphi(0)+ \psi(0) \bigr) \end{pmatrix}}, \
 {\begin{pmatrix} \theta_3(0) & 0 \\ 0 & \theta_3(0) \end{pmatrix}}, \dots, {\begin{pmatrix} \theta_l(0) & 0 \\ 0 & \theta_l(0) \end{pmatrix}} \right) \\
& = \left( {\begin{pmatrix} 0 & 1 \\ -\varphi(0)\psi(0) & t \bigl( 1 + \varphi(0)\psi(0) \bigr) \end{pmatrix}}, \
 {\begin{pmatrix} 0 & 1 \\ -\varphi(0)\psi(0) & t \bigl( 1 + \varphi(0)\psi(0) \bigr) \end{pmatrix}}, \
 {\begin{pmatrix} \theta_3(0) & 0 \\ 0 & \theta_3(0) \end{pmatrix}}, \dots, {\begin{pmatrix} \theta_l(0) & 0 \\ 0 & \theta_l(0) \end{pmatrix}} \right)
\end{align*}
to show that the right-hand side of an equality similar to the one as in the claim in the proof of Theorem \ref{multilinear} vanishes.
This proves (\ref{skewsymmetryeq}).

From (\ref{skewsymmetryeq}), we have, using (\ref{linear-equa}) and (\ref{linear-equb}),
$$ \bigl( \varphi(t)\psi(t),\, \varphi(t)\psi(t), \theta_3(t), \dots, \theta(l) \bigr)
= \bigl(\varphi(t),\, \varphi(t), \theta_3(t), \dots, \theta(l) \bigr) + \bigl( \psi(t),\, \psi(t), \theta_3(t), \dots, \theta(l) \bigr).$$

On the other hand, by Theorem \ref{multilinear}, we also have
\begin{multline*}
\quad \bigl( \varphi(t)\psi(t),\, \varphi(t)\psi(t), \theta_3(t), \dots, \theta(l) \bigr) \\
= \bigl( \varphi(t),\, \varphi(t),\theta_3(t), \dots, \theta(l) \bigr) + \bigl( \varphi(t),\, \psi(t),\theta_3(t), \dots, \theta(l) \bigr ) \\
                + \bigl ( \psi(t), \, \varphi(t),\theta_3(t), \dots, \theta(l) \bigr) + \bigl( \psi(t), \, \psi(t),\theta_3(t), \dots, \theta(l) \bigr).
\end{multline*}

The equality of the right hand sides of these two identities leads to the skew-symmetry.
\end{proof}

The irreducibility assumption in Proposition \ref{skewsymmetry} can be replaced by an assumption similar to one of the conditions
in Corollary \ref{multilinear1} or \ref{multilinear2}. For example, it is enough to require that the symbol $\bigl( \theta_1(t),\dots, \theta_l(t) \bigr)$
is irreducible if the field $k$ is of characteristic 0.

     \bibliographystyle{plain}

\end{document}